\documentclass[a4paper,12pt]{amsart}
\usepackage[mathscr]{eucal}
\usepackage{fullpage}

\def\N{\mathbb{N}}
\def\E{{\bf E}}
\def\PP{\mathbb{P}}

\def\F{\mathcal{F}}

\def\a{\alpha}

\def\k{\kappa}

\def\o{\omega } \def\O{\Omega }

\def\i{\infty}
\def\mylim{\lim_{n\to\infty}}

\newcommand{\f}{\varphi}
\newcommand{\e}{\varepsilon}
\newcommand{\Reals}{{\mathbb{R}^1}}
\newcommand{\condex}[1]{\E\left[\left.#1\right|\F_i\right]}
\newcommand{\condexn}[1]{\E\left[\left.#1\right|\F_n\right]}

\newtheorem{theorem}{Theorem}
\newtheorem{lemma}[theorem]{Lemma}

\newtheorem{corollary}[theorem]{Corollary}
\newtheorem{assumption}{Assumption}

\theoremstyle{definition}

\theoremstyle{remark}
\newtheorem{remark}[theorem]{Remark}

\newcommand{\wt}{\widetilde}

\title[decay rates in stochastic
  difference equation]
{Non-exponential stability and decay rates in nonlinear stochastic
  difference equation with unbounded noises}

\thanks{The first author is supported by an Albert College Fellowship, awarded by the Dublin
City University Research Advisory Panel. The third author is
supported by the Mona Research Fellowship Programm awarded by the
University of the West Indies, Mona Campus, Jamaica.}
\author[J.~Appleby]{John A. D. Appleby}
\address{School of Mathematical Sciences, Dublin City
University, Dublin 9, Ireland} \email{john.appleby@dcu.ie}
\urladdr{http://webpages.dcu.ie/\textasciitilde applebyj}
\author[G.~Berkolaiko]{Gregory Berkolaiko}
\address{Department of Mathematics, Texas A\&M University, USA}
\email{gregory.berkolaiko@math.tamu.edu}
\urladdr{http://www.math.tam.edu/\textasciitilde berko}
\author[A.~Rodkina]{Alexandra Rodkina}
\address{Department of Maths\&CSci,The
University of the West Indies, Kingston, Jamaica. }
\email{alexandra.rodkina@uwimona.edu.jm}

\begin{document}
\begin{abstract}
We consider stochastic difference equation
  \begin{displaymath}
  x_{n+1}=x_n\biggl(1-hf(x_n)+\sqrt{h} g(x_n) \xi_{n+1}\biggr), \quad n=0, 1,
  \dots,\quad x_0\in \Reals,
  \end{displaymath}
  where functions $f$ and $g$ are nonlinear and bounded, random variables
  $\xi_{i}$ are independent and $h>0$ is a nonrandom parameter.

  We establish results on asymptotic stability and instability of the trivial
  solution $x_n\equiv0$.  We also show, that for some natural choices of the
  nonlinearities $f$ and $g$, the rate of decay of $x_n$ is approximately
  polynomial: we find $\alpha>0$ such that $x_n$ decays faster than
  $n^{-\alpha+\varepsilon}$ but slower than $n^{-\alpha-\varepsilon}$ for any
  $\varepsilon>0$.

  It also turns out that if $g(x)$ decays faster than $f(x)$ as $x\to0$, the
  polynomial rate of decay can be established exactly, $x_n n^\alpha \to
  \textrm{const}$.  On the other hand, if the coefficient by the noise does not decay
  fast enough, the approximate decay rate is the best possible result.
\end{abstract}
\maketitle

%%%%%%%%%%%%%%%%%%%%%%%%%%%%%%%%%%%%%%%%%%%%%%%%%%%%%%%%%%%
%
%
%\documentclass[12pt]{article}
%
%\textwidth 160mm \textheight 220mm \oddsidemargin 0mm \topmargin
%-10mm
%
%\usepackage{amsmath,amsthm,amssymb}
%\usepackage[mathscr]{eucal}
%\usepackage{amssymb}
%\usepackage[notref,notcite]{showkeys}
%
%\usepackage{amsthm}
%\newtheorem{theorem}{Theorem}
%\newtheorem{lemma}{Lemma}
%\newtheorem{cor}{Corollary}

%%%%%%%%%%%%%%%%%%%%%%% SECTION %%%%%%%%%%%%%%%%%%%%%%%%%%%%%%%%%
\section{Introduction}

In this paper we address the questions of stability and the rate
of decay of solutions of the difference equation
\begin{equation}
  \label{eq:our_eq}
  x_{n+1} = x_n \biggl(1-hf(x_n)+\sqrt{h} g(x_n) \xi_{n+1}\biggr),
  \quad n=0,1,\dots,\quad x_0\in \Reals,
\end{equation}
where $\xi_{n+1}$ are independent random variables.  The functions
$f$ and $g$ are nonlinear and are assumed to be bounded.  The
small parameter $h>0$ usually arises as the step size in numerical
schemes. Equation (\ref{eq:our_eq}) may be viewed as
stochastically perturbed version of a deterministic autonomous
difference equation, where the random perturbation is
state-dependent.  In general, it does not have linear leading
order spatial dependence close to the equilibrium.  As a
consequence of the non-hyperbolicity of the equilibrium, the
convergence of solutions of (\ref{eq:our_eq}) to its equilibrium
zero cannot be expected to take place at an exponentially fast
rate.

Similarly to deterministic difference equations, analyzing
asymptotic behavior of stochastic difference equations is often
harder (see \cite{AMR, AMaR, BR,
  KS, Rodkina01,RM01, RS, RS1}) than analyzing their {\em differential\/}
counterparts.  Nonetheless, we feel it is very important to
develop techniques and better understanding of the similarities
and the differences between the two types of equations.  In this,
we are motivated by two principal reasons. Firstly, in many
applied contexts the studied phenomena are intrinsically discrete
(see, for example, \cite{Saber}).  Using continuous approximation
can sometimes mask significant phenomena.  Going in the other
direction, numerical simulation of stochastic differential
equations involves solving an associated difference equation.  It
is important to know whether the discretization can produce
spurious behaviors and how can this be avoided.  For example, one
needs to study if the asymptotic behavior of the discretized
equation is a faithful reproduction of the asymptotic behavior of
the original equation. The corresponding property is called
``A-stability'' and it has been addressed in the stochastic
context in \cite{Hig00}, \cite{KP}, \cite{SM96}.

In this paper we analyze sufficient and necessary conditions for
solutions $x_n$ to converge to zero as $n\to\infty$
(``stability'') and the rate at which such convergence happens for
different types of the nonlinearities $f$ and $g$.  Our results
should be compared to an earlier work \cite{ARS} (see also
\cite{AR}) in which similar differential equations had been
analyzed.

One of the technical difficulties arising in the study of
stability of stochastic difference equations is dealing with
unbounded noise.  So far, many results have been only available
for the case of bounded noises (e.g.~\cite{BR}).  Yet, one of the
most applicable scenarios, discretization of the white noise,
involves normally distributed (and thus unbounded) random
variables.  In this paper we develop a tool which is designed to
overcome this difficulty.  In particular, it is instrumental in
proving the instability result (Theorem~\ref{thm:nonlinunstab} in
Section~\ref{sec:stab}) in this paper and can also be used to
prove instability in several related models (e.g.\ in
\cite{HMS,HMY}).  Section~\ref{sec:ito_form} is devoted to setting
up this tool, which can be thought of as a discrete variant of the
It\=o formula, and proving the corresponding theorem
(Theorem~\ref{thm:my_ito}).

Armed with Theorem~\ref{thm:my_ito} we formulate and prove
criteria for almost sure asymptotic stability and instability of
solutions to equation~(\ref{eq:our_eq}).  In
Section~\ref{sec:decay_rate} we concentrate on the decay rate of
the solutions (assuming they converge to 0).  The principal result
here is the comparison theorem which provides implicit information
on asymptotic behavior of solution $x_n$ via the limit
\begin{equation*}
  \lim_{n\to\infty} \frac{\ln|x_n|}{\sum_{i=1}^n S(x_i)},
\end{equation*}
where $S(x)$ stands for either $g^2(x)$ or $|f(x)|$.  In the
special (but typical) cases of polynomially decaying $f$ and $g$,
we extract explicit information (see Corollary~\ref{cor:rate}) on
the decay rate of $x_n$ in the form of the limit
\begin{equation*}
  \lim_{n\to\infty} \frac{\ln|x_n|}{\ln n} = -\lambda < 0.
\end{equation*}
The above limit allows one to conclude that the decay rate of
$x_n$ is of polynomial type.  More precisely, for any
$\epsilon>0$, the following bound is valid {\em eventually\/} as
$n\to\infty$,
\begin{equation}
  \label{eq:poly_bound}
  n^{-\lambda-\epsilon} \leq |x_n| \leq n^{-\lambda+\epsilon}.
\end{equation}

At this point a natural question arises: under what circumstances
bound~(\ref{eq:poly_bound}) can be strengthened to the exact
power-law decay $x_nn^\lambda \to \textrm{const}$? This question
is answered in Section~\ref{sec:exact_rate}. Heuristically, the
answer can be described as follows. The convergence to zero can be
caused either by the deterministic term $f(x)$ or by the noise
term $g(x)\xi$, depending on the comparative speed of decay of
$f(x)$ and $g(x)$ as $x$ tends to 0. When $f(x)$ is dominant, the
convergence of $x_n$ happens at an exact power-law rate. On the
other hand, if the noise term is significant, we show that an
exact rate result is {\em
  impossible}, namely
\begin{equation*}
  \limsup_{n\to \infty}|x_n|n^\lambda = \infty,
  \qquad\mbox{and}\qquad
  \liminf_{n\to \infty}|x_n|n^\lambda = 0,
\end{equation*}
for some $\lambda$.

%%%%%%%%%%%%%%%%%%%%%%% SECTION %%%%%%%%%%%%%%%%%%%%%%%%%%%%%%%%%
\section{Auxiliary Definitions and Facts}
In this section we give a number of necessary definitions and a
lemmas we use to prove our results. A detailed exposition of the
definitions and facts of the theory of random processes can be
found in, for example, \cite{Shiryaev96}.

Let $(\Omega, {\mathcal{F}}, \{{\mathcal{F}}_n\}_{n \in \N},
{\PP})$ be a complete filtered probability space. Let
$\{\xi_{n}\}_{n\in\N}$ be a sequence of independent random
variables with $\E \xi_{n}=0$. We assume that the filtration
$\{{\mathcal{F}}_n\}_{n
  \in \N}$ is naturally generated: $\mathcal{F}_{n+1} = \sigma \{\xi_{i+1} :
i=0,1,...,n\}$.

Among all the sequences $\{X_n\}_{n \in \N}$ of the random
variables we distinguish those for which $X_n$ are
${\mathcal{F}}_n$-measurable $\forall n \in \N$.

A stochastic sequence $\{X_n\}_{n \in \N}$ is said to be an {\em
  $\mathcal{F}_n$-martingale}, if ${\bf E}|X_n|<\i$ and $\E\bigl[X_n\bigl|
\mathcal{F}_{n-1}\bigr]=X_{n-1}$ for all $n\in\N$ {\it a.s.}

A stochastic sequence $\{\xi_n\}_{n \in N}$ is said to be an {\em
$\mathcal{F}_n$-martingale-difference}, if $\E|\xi_n|<\i$ and
$\E\bigl(\xi_n\bigl| \mathcal{F}_{n-1}\bigr)=0$ {\it a.s.} for all
$n\in\N$.

We use the standard abbreviation ``{\it a.s.}'' for the wordings
``almost sure'' or ``almost surely'' throughout the text.

%\begin{example}
%  Let $\{X_n\}_{n\in \N}$ be a sequence of $\mathcal{F}_n$-measurable random variables
%  such that $\E\bigl(X_n\bigl|\mathcal{F}_{n-1}\bigr)= 1$. Let $Z_n=\prod_{i=1}^nX_i$
%  and $\E |Z_n|<\i$ for all $n\in\N$.  Then $\{Z_n\}_{n\in \N}$ is a
%  martingale.
%\end{example}

If $\{X_n\}_{n\in \N}$ is a martingale, in the form
$X_n=\sum_{i=1}^n \rho_i$, then the \emph{quadratic variation} of
$X$ is the process $\langle X\rangle$ defined by
\[
\langle X_n\rangle=\sum_{i=1}^n
\mathbf{E}[\rho_i^2|\mathcal{F}_{i-1}].
\]

Three lemmas below are variants of martingale convergence theorems
(see e.g. \cite{Shiryaev96}).

\begin{lemma}\label{lem:conv_variance}
  If $\{X_n\}_{n\in \N}$ is a martingale, $X_n=\sum_{i=1}^n \rho_i$, then
  \begin{displaymath}
    \left\{\omega: \sum_{i=1}^\infty \mathbf{E}[\rho_i^2|\mathcal{F}_{i-1}] < \i \right\}
    \subseteq \{X_n \to \}.
  \end{displaymath}
  Here $\left\{X_n \to \right\}$ denotes the set of all
  $\omega\in\Omega$ for which $\lim\limits_{n\to\infty}X_n$ exists
  and is finite.
\end{lemma}

\begin{lemma}\label{lem:div_variance}
  If $\{X_n\}_{n\in \N}$ is a martingale, $X_n=\sum_{i=1}^n \rho_i$, and
  \begin{displaymath}
    \sum_{i=1}^\infty \mathbf{E}[\rho_i^2\big|\mathcal{F}_{i-1}] = \i, \quad a.s.
  \end{displaymath}
  Then, a.s.,
  \begin{displaymath}
    \frac{X_n}{\sum_{i=1}^n \mathbf{E}\left[\rho_i^2\big|\mathcal{F}_{i-1}\right]} \to 0, \qquad n\to\infty.
  \end{displaymath}
\end{lemma}

\begin{lemma}\label{lem:pos_submart}
%  If $\{X_n\}_{n\in \N}$ is a non-negative submartingale, $X_n=\sum_{i=1}^n
%  \rho_n$, then
%  \begin{displaymath}
%    \left\{\sum_{i=1}^\infty \condex{\rho_i} < \i \right\}
%    \subseteq \{X_n \to \}.
%  \end{displaymath}
%  In particular,
  If $X_n$ is a non-negative martingale, then $\mylim X_n$ exists with
  probability 1.
\end{lemma}

The following lemma is proved in \cite{AMaR}.
\begin{lemma}\label{lem:nonegdif}
  Let $\{Z_n\}_{n\in \N}$ be a non-negative  $\mathcal{F}_n$-measurable
  process, ${\E}|Z_n|<\i$  $\forall n\in \N$, and
  \begin{displaymath}
    Z_{n+1}\le Z_{n}+u_n-v_n+\nu _{n+1}, \quad n = 0, 1, 2, \dots,
  \end{displaymath}
  where $\{\nu_n\}_{n\in \N}$ is an ${\mathcal{F}}_n$-martingale-difference,
  $\{u_n\}_{n\in \N}$, $\{v_n\}_{n\in \N}$ are nonnegative $\mathcal{F}_n$-measurable
  processes and ${\E}|u_n|, {\E}|v_n|<\i$ $\forall n\in \N$.

  Then $$\left\{\o: \sum_{n=1}^{\i} u_n<\i\right\}\subseteq \left\{\o:
    \sum_{n=1}^{\i} v_n<\i\right\}\bigcap\{Z\to\}.$$
\end{lemma}

%%%%%%%%%%%%%%%%%%%%%%% SECTION %%%%%%%%%%%%%%%%%%%%%%%%%%%%%%%%%
\section{Discretized It\={o} formula}
\label{sec:ito_form}

Below we will make use of the notation $O(\cdot)$:
\begin{equation}
  \label{eq:big_o}
  \alpha(r) = O(\beta(r)) \quad \mbox{as}\quad r\to r_0
  \qquad \Leftrightarrow \qquad
  \limsup_{r\to r_0} \left|\frac{\alpha(r)}{\beta(r)}\right| < \infty.
\end{equation}
Here $r_0$ can be a real number or $\pm\infty$ and the argument
$r$ can be both continuous and discrete.  We will also use
$o(\cdot)$:
\begin{equation*}
  \alpha(r) = o(\beta(r)) \quad \mbox{as}\quad r\to r_0
  \qquad \Leftrightarrow \qquad
  \lim_{r\to r_0} \frac{\alpha(r)}{\beta(r)} = 0.
\end{equation*}

\begin{assumption}
  \label{assum:noise}
  We will make the following assumptions about the noise $\xi_n$:
  \begin{enumerate}
  \item $\xi_n$ are independent random variables satisfying
    \begin{equation*}
      \E\xi_n = 0, \qquad \E\xi_n^2=1, \qquad
      \E|\xi_n|^3 \mbox{ are uniformly bounded,}
    \end{equation*}
  \item the probability density functions $p_n(\xi)$ exist and satisfy
    \begin{equation*}
      x^3 p_n(x) \to 0 \quad \mbox{as } |x|\to\infty \quad
      \mbox{uniformly in } n.
    \end{equation*}
  \end{enumerate}
\end{assumption}

The following theorem can be thought of as a discretized relative
of the It\=o formula.

\begin{theorem}
  \label{thm:my_ito}
  Consider $\f:\Reals\to\Reals$ such that there exists $\delta>0$ and
  $\wt{\f}:\Reals\to\Reals$ satisfying
  \begin{enumerate}
  \item $\wt{\f} \equiv \f$ on $U_\delta = [1-\delta,1+\delta]$,
  \item $\wt{\f}\in C^3(\Reals)$ and $|\wt{\f}'''(x)|\leq M$ for some $M$
    and all $x\in\Reals$,
  \item $\int_{\Reals} |\f - \wt{\f}| dx <
    \infty$.
  \end{enumerate}
  Let $f$ and $g$ be $\mathcal{F}$-measurable bounded random variables; $\xi$ be an
  $\mathcal{F}$-independent random variable satisfying Assumption~\ref{assum:noise}.
  Then
  \begin{equation}
    \label{eq:ito_form}
    \E \left[\left. \f\left(1+fh + g\sqrt{h}\xi\right) \right| \mathcal{F} \right ]
    = \f(1) + \f'(1) f h + \frac{\f''(1)}{2} g^2 h
    + hf o(1) + hg^2 o(1),
  \end{equation}
  where the error terms $o(1)$ satisfy
  \begin{enumerate}
  \item if $|f|,|g| < K$ then $o(1) \to 0$ as $h\to0$, uniformly in $f$ and $g$,
  \item if $h<H$ then $o(1) \to 0$ as $f\to0$ and $g\to0$ uniformly in $h$.
  \end{enumerate}
\end{theorem}

\begin{proof}
  For brevity we will assume that $f$ and $g$ are constants and
  correspondingly use the non-conditional expectation, the proof of the
  general case being completely analogous.

  The proof consists of two main parts.  In the first part we derive
  formula~(\ref{eq:ito_form}) for $\E\left[\wt{\f}\right]$.  In the second
  part we prove that, for our purposes, $\wt{\f}$ is a good approximation for
  $\f$.  More precisely, we prove the following estimate for the error term,
  \begin{equation*}
    \E\left[\f - \wt{\f}\right]
    = h g^2 o(1).
  \end{equation*}

  Part A. By Taylor expansion,
  \begin{equation*}
    \wt{\f}(1+x) = \wt{\f}(1) + \wt{\f}'(1)x + \frac{\wt{\f}''(1)}{2}x^2
    + \frac{\wt{\f}'''(\theta)}{6}x^3,
  \end{equation*}
  with $\theta$ lying between $0$ and $x$.  We substitute $x = fh +
  g\sqrt{h}\xi$ and take expectation.  Using properties of $\xi$,
  \begin{equation*}
    \E x = fh, \qquad \E x^2 = f^2h^2 + g^2h,
  \end{equation*}
  and therefore
  \begin{equation*}
    \E\wt{\f}(x) = \wt{\f}(1) + \wt{\f}'(1) f h + \frac{\wt{\f}''(1)}{2} g^2 h
    + \wt{\f}''(1)f^2h^2/2 + \E[ \wt{\f}'''(\theta) x^3] / 6.
  \end{equation*}
  Because $\wt{\f}'''(\theta)$ is uniformly bounded, we can estimate, by
  expanding $x^3$,
  \begin{equation*}
    \left|\E[ \wt{\f}'''(\theta) x^3] / 6\right|
    \leq M \E|x^3|/ 6 \leq f h^2 O(f h^{1/2}) + g^2 h O(g h^{1/2})+g^2hO(fh).
  \end{equation*}
  This proves formula (\ref{eq:ito_form}) for $\wt{\f}\left(1+fh +
    g\sqrt{h}\xi \right)$.

  Part B.  We introduce the shorthand $c_1=1+hf$ and $c_2=\sqrt{h}g$ and seek
  an estimate for the error term
  \begin{equation*}
    \Delta = \E\left[ \f\left(c_1+c_2\xi\right)
      - \wt{\f}\left(c_1+c_2\xi\right)\right].
  \end{equation*}
  We have
  \begin{eqnarray*}
    \Delta &=& \int_{\Reals} \left(\f\left(c_1+c_2\xi\right)
      - \wt{\f}\left(c_1+c_2\xi\right)\right) p(\xi) d\xi \\
    &=& \int_{\Reals} \left(\f(r)-\wt{\f}(r)\right)
    p\left(\frac{r-c_1}{c_2}\right) \frac{dr}{|c_2|}
    = \int_{\Reals\setminus U_\delta} \left(\f(r)-\wt{\f}(r)\right)
    p\left(\frac{r-c_1}{c_2}\right) \frac{dr}{|c_2|},
  \end{eqnarray*}
  where we introduced a new variable of integration, $r = c_1+c_2\xi$, and
  excluded $U_\delta$ from the integration range because $\f(r)-\wt{\f}(r) =
  0$ on $U_\delta$.

  Now we can estimate
  \begin{eqnarray*}
    |\Delta| &\leq& \sup_{r\not\in U_\delta}
    \left\{p\left(\frac{r-c_1}{c_2}\right) \frac{1}{|c_2|}\right\}
    \int_\Reals \left|\f(r)-\wt{\f}(r)\right|dr\\
    &\leq& |c_2|^2 C \sup_{r\not\in U_\delta}
    \left\{p\left(\frac{r-c_1}{c_2}\right) \frac{1}{|c_2|^3}\right\}
    = h g^2 C \sup_{r\not\in U_\delta}
    \left\{\frac{p(y)y^3}{(r-1-hf)^3}\right\},
  \end{eqnarray*}
  where
  \begin{equation*}
    y = \frac{r-1-hf}{\sqrt{h}g}.
  \end{equation*}
  If {\em either\/} $h$ is bounded and $f,g\to0$ {\em or\/} $|f|$ and $|g|$
  are bounded and $h\to0$, it is easy to see that $y\to\infty$ uniformly on
  $r\in \Reals\setminus U_\delta$.  Since under the same conditions
  $(r-1-hf)^3$ is bounded away from zero, the assumption $p(y)y^3 \to 0$
  implies
  \begin{equation*}
    \sup_{r\not\in U_\delta} \left\{\frac{p(y)y^3}{(r-1-hf)^3}\right\}
    = o(1),
  \end{equation*}
  and, therefore, $|\Delta| \leq h g^2 o(1)$.
\end{proof}

%%%%%%%%%%%%%%%%%%%%%%%%%%%%% SECTION %%%%%%%%%%%%%%%%%%%%%%%%%%%%%%%%%%%%%
\section{Stability and instability}
\label{sec:stab}

We consider equation
\begin{eqnarray}
  x_{n+1}=x_n\left(1+h f(x_n)+\sqrt{h} g(x_n)\xi_{n+1}\right), \quad
  n=0,1, \dots,  \label{eq0}
\end{eqnarray}
with nonrandom initial value $x_0\in \Reals$, and independent
random variables $\xi_n$ satisfying $\E\xi_n=0$, $\E\xi_n^2=1$ for
all $n\in \N$. The functions $g, f: \Reals\to \Reals$ are
nonrandom, continuous and bounded:
\begin{equation}
  \label{eq:boundfg}
  |g(u)|, |f(u)|\le 1\quad \forall u\in\Reals.
\end{equation}

\begin{theorem}\label{thm:stability}
  Let functions $f$ and $g$ be bounded and $\xi_n$ satisfy
  Assumption~\ref{assum:noise}.  Let also
  \begin{equation}
    \sup_{u\in \Reals \backslash \emptyset}\left\{\frac{2f(u)}{g^2(u)}
    \right\} = \beta < 1.
    \label{eq:ratio_fg}
  \end{equation}
  If $h$ is small enough then $\mylim x_n(\o)=0 $ {\it a.s.} where $x_n$ is a
  solution to equation (\ref{eq0}).
\end{theorem}

\begin{remark}
  If $g(u)=0$ for some $u\neq0$, we consider (\ref{eq:ratio_fg})
  fulfilled iff $f(u)<0$.  Thus we impose no restrictions on $g(u)$ when
  $f(u)<0$ for all nonzero $u$.
\end{remark}

%%%%%% Proof Stability
\begin{proof}
  We raise equation~(\ref{eq0}) to a power $\alpha>0$, which will be
  determined later,
  \begin{equation*}
    |x_{n+1}|^\alpha
    = |x_n|^\alpha \left|1 + h f(x_n) + \sqrt{h} g(x_n)\xi_{n+1}\right|^\alpha,
  \end{equation*}
  and denote $z_n = |x_n|^\alpha$.  We define $\phi_\alpha(y)=|y|^\alpha$,
% by
%  \begin{equation}
%    \label{eq:def_phi_alpha}
%    \phi_\alpha \in C^3(\Reals), \qquad
%    \phi_\alpha(y) =
%    \begin{cases}
%      |1+y|^\alpha, &y\not\in (-1.5, 0.5),\\
%      \phi_\alpha(y)\geq |1+y|^\alpha,  &y\in (-1.5, 0.5).
%    \end{cases}
%  \end{equation}
%  It can be easily seem that $\phi_\alpha'''$ is bounded on $\Reals$.  Then
%  \begin{equation*}
%    z_{n+1} \leq z_n\: \phi_\alpha \big(h f(x_n)
%      + \sqrt{h} g(x_n)\xi_{n+1}\big).
%  \end{equation*}
  denote
  \begin{eqnarray}
    \label{eq:def_Phi}
    \Phi_n &=& \condexn{\phi_\alpha \big(1 + h f(x_n)
      + \sqrt{h}g(x_n)\xi_{n+1}\big)} - 1,\\
    \label{eq:def_rho}
    \rho_{n+1} &=& z_n \left(\phi_\alpha \big(1 + h f(x_n)
      + \sqrt{h}g(x_n)\xi_{n+1}\big) - \Phi_n - 1\right)
  \end{eqnarray}
  and rewrite
  \begin{equation}
    \label{eq:x_alpha_est}
    z_{n+1} = z_n + z_n \Phi_n + \rho_{n+1}.
  \end{equation}
  Here $\rho_{n+1}$ is an $\mathcal{F}_{n+1}$-martingale-difference.  Applying
  Theorem~\ref{thm:my_ito} to $\Phi_n$ we have: $\phi_\alpha(1)=1$,
  $\phi_\alpha'(1)=\alpha$, $\phi_\alpha''(1)=\alpha(\alpha-1)$ and, therefore,
  \begin{eqnarray}
    \label{eq:Phi_expr1}
    \Phi_n &=&  \alpha h f(x_n)\big(1+o(1)\big)
    + \frac{\alpha(\alpha-1)}{2} h g^2(x_n) \big(1+o(1)\big)\\
    \label{eq:Phi_expr2}
    &\leq& \frac12 \alpha h g^2(x_n)
    \left(\beta + \alpha - 1 + o(1)\right),
  \end{eqnarray}
  where we used condition~(\ref{eq:ratio_fg}) to obtain the estimate on the
  second line.  Since $\beta < 1$, we can choose $\alpha$ and $h$ sufficiently
  small, so that $\Phi_n \leq 0$ for all $x_n$.  Now we can apply
  Lemma~\ref{lem:nonegdif} with $u_n=0$ and $v_n=-z_n \Phi_n$ and conclude
  that $z_n$ converges to some (possibly random) value $z_\infty$.

  To prove that $z_\infty$ is a.s.~0 we assume the contrary: there exists
  $y\neq0$ such that, for any $\delta>0$, the probability that the limit of
  $z_n$ lies in the interval $(y-\delta,y+\delta)$ is nonzero.

  At the point $y$ we either have $g(y)\neq0$ and then $\Phi_n<0$ by
  (\ref{eq:Phi_expr2}) or $g(y)=0$ but then $f(y)<0$ and again $\Phi_n<0$ (now
  using (\ref{eq:Phi_expr1})).  By continuity, $\Phi_n$ remains bounded away
  from zero in a $\delta$-neighborhood of $y$.  Therefore $\sum_{n=1}^\infty
  z_n\Phi_n$ is divergent which contradicts Lemma~\ref{lem:nonegdif}.
\end{proof}

It turns out that condition (\ref{eq:ratio_fg}) is close to being
necessary for stability.  For the same equation we now ask the
opposite question: under what conditions on $f$ and $g$ solutions
of (\ref{eq0}) do {\it not\/} tend to zero.

%%% Instability Theorem
\begin{theorem}\label{thm:nonlinunstab}
  Let $f$ and $g$ be bounded and $\xi_n$ satisfy
  Assumption~\ref{assum:noise}. Let also
  \begin{displaymath}
    f(u)>0 \quad \mbox{and} \quad g(u)\neq0 \quad \mbox{when}\quad  u\neq 0
  \end{displaymath}
  and
  \begin{equation}
    \label{eq:limneq0}
    \liminf_{u\to 0}\left\{\frac{2|f(u)|}{g^2(u)}\right\} > 1.
  \end{equation}
  If $x_n$ be a solution to equation (\ref{eq0}) with an initial value
  $x_0\in\Reals$ and $h$ is small enough then
  $\PP\left\{\mylim x_n(\o) = 0\right\} = 0$.
\end{theorem}

%%% Proof of Instability Theorem
\begin{proof}
Consider
\begin{equation}
  \label{eq:cond_exp}
  \Phi_i = \E\left[\left|1 + hf(x_i) + \sqrt{h}g(x_i)\xi_{i+1}\right|^{-\a}
    \biggl|\mathcal{F}_i\right]
\end{equation}
with $\alpha<1$ and with $\mathcal{F}_0$ being the trivial
$\sigma$-algebra. Since $\f(x) = |x|^{-\alpha}$ is integrable
around $x=0$ and has bounded third derivative outside a
neighborhood of $0$, we can apply Theorem~\ref{thm:my_ito} to
obtain
\begin{equation*}
  \Phi_i = 1 -\a h f(x_n)\big(1+o(1)\big)
    + \frac{\a(\a+1)}{2} h g^2(x_n)\big(1+o(1)\big).
\end{equation*}
In particular, $\Phi_i$ are finite.  Therefore, we can form
\begin{equation*}
  M_n = \prod_{i=0}^{n-1}
  \frac{\left|1-hf(x_i)+\sqrt{h}g(x_i)\xi_{i+1}\right|^{-\a}}{\Phi_i}.
\end{equation*}
 Since
$\E|M_n|=\E M_n = 1<\infty$, $M_n$ is a positive martingale,
convergent by Lemma~\ref{lem:pos_submart}.

We now have the following representation for the solution $x_n$
\begin{equation}
  \label{eq:prod_rep}
  |x_{n}|^\a
  = |x_0|^\a \prod_{i=0}^{n-1}
  \left|1 - hf(x_i)+\sqrt{h}g(x_i)\xi_{i+1}\right|^\a
  = |x_0|^\a M_n^{-1} \prod_{i=0}^{n-1} \Phi_i^{-1}.
\end{equation}

Suppose now that the theorem is untrue, i.e.~that $\PP(\O_1)>0$,
where $\O_1=\{\o: \lim x_n(\o)=0\}$.  Using condition
(\ref{eq:limneq0}), for each $\o\in\O_1$ we can find $\delta\in
(0,1)$ and $N(\o, \delta)$ such that $g^2(x_n) <
2f(x_n)(1-\delta)$ for all $n>N(\o,\delta)$.

Thus, for $n>N(\o,\delta)$,
\begin{equation*}
  \Phi_i < 1 + \a hf(x_n)\Big( -1 + o(1)
    + (\a+1)(1-\delta)\big(1+o(1)\big)\big) < 1,
\end{equation*}
when $h$ and $\a$ are small enough.

Applying this inequality to representation~(\ref{eq:prod_rep}) we
obtain
\begin{equation*}
  |x_{n}|^\a
  \geq |x_0|^\a M_n^{-1} \prod_{i=0}^{N(\o,\delta)} \Phi_i^{-1}.
\end{equation*}
The only $n$-dependent factor on the right-hand side is $M_n^{-1}$
which tends to a nonzero limit by Lemma~\ref{lem:pos_submart}. All
other factors being nonzero as well, we conclude that $x_n^\a$
remains bounded away from 0, which contradicts our definition of
$\O_1$.
\end{proof}

\begin{remark}
  A more refined analysis of the factors in representation~(\ref{eq:prod_rep})
  allows to strengthen the conclusion of Theorem~\ref{thm:nonlinunstab} to
  $\PP\left\{\liminf_{n\to \i}x_n(\o)=0\right\} = 0$.  The proof of this
  statement, however, is unpleasantly technical and we leave it out.
\end{remark}

%%%%%%%%%%%%%%%%%%%%%%%%%%%%%%% SECTION %%%%%%%%%%%%%%%%%%%%%%%%%%%%%%%%%%%
\section{Decay Rate}
\label{sec:decay_rate}

In this section we establish results on the a.s.\ decay rate of
solutions $x_n$ of (\ref{eq0}).

The first subsection contains some variations of the classical
Toeplitz lemma from analysis.  In the second subsection we present
a result about asymptotic behavior of $\ln x^2_n$ in two cases:
when $\lim_{u\to
  0}\frac{f(u)}{g(u)^2}=L$ and when $\lim_{u\to 0}\frac{f(u)}{g(u)^2}=-\i$.
These conditions include, but are weaker than, the sufficient
conditions given by Theorem~\ref{thm:stability} for the stability
of $x_n$.  For this reason we explicitly assume that $x_n\to0$
when we discuss the rate of decay.

%%%%%%%%%%%%%%%%%%%%%%%%%%%% Subsection %%%%%%%%%%%%%%%%%%%%%%%%%%%%%%%%%%
\subsection{Variations on Toeplitz Lemma}

In this section we state Toeplitz Lemma and prove one of its
corollaries. The version of Toeplitz Lemma we need is taken from
(\cite{Shiryaev96}, p. 390).

\begin{lemma}[Toeplitz Lemma]\label{lem:toeplitz}
  Let $\{a_n\}_{n\in\N}$ be a sequence of nonnegative real numbers such that
  $\sum_{i=1}^\infty a_i$ diverges.  If $\k_n\to\k_\i$ as $n \to \i$ then
  \begin{equation*}
    \mylim \frac{\sum_{i=0}^n a_i\k_i}{\sum_{i=0}^n a_i} = \k_\i .
  \end{equation*}
\end{lemma}

We will use the following 2 corollaries of Toeplitz Lemma.

\begin{lemma}\label{lem:toepl_cor1}
  Let $\{a_n\}_{n\in \mathbb{N}}$ be a sequence of nonnegative real numbers with
  $\sum_{i=1}^na_i\to \i$ when $n\to\i$.  Then
  \begin{displaymath}
    \mylim \frac{b_n}{a_n}=c \quad \Rightarrow \quad
    \mylim \frac{\sum_{i=1}^nb_i}{\sum_{i=1}^na_i}=c.
  \end{displaymath}
  Also,
\begin{displaymath}
    \limsup_{n\to
    \infty}\frac{\sum_{i=1}^nb_i}{\sum_{i=1}^na_i}=\infty
    \quad \Rightarrow \quad \limsup_{n\to \infty} \frac{b_n}{a_n}=\infty.
  \end{displaymath}
\end{lemma}
\begin{proof}
  The statement follows from the representation
  \begin{displaymath}
    \mylim \frac{\sum_{i=1}^nb_i}{\sum_{i=1}^na_i}
    = \mylim \frac{\sum_{i=1}^n\frac{b_i}{a_i}a_i}{\sum_{i=1}^na_i},
  \end{displaymath}
  and Toeplitz Lemma.
\end{proof}

The following lemma is useful for extracting information about a
sequence $\{y_n\}_{n\in \mathbb{N}}$ if what is known is given in
an implicit form such as \, $f(y_n)(y_{n+1}-y_n)\to c$ \, for some
function $f$.
\begin{lemma}\label{lem:T}
  Let $c>0$, let $f:\Reals_+\to \Reals_+$ be monotonous continuous
  function. Let $\{y_n\}_{n\in \mathbf{N}}$ be a positive increasing sequence
  such that $ \mylim \frac{f(y_n)}{f(y_{n-1})}=1$ and $\mylim \Delta y_n =0$,
  where $\Delta y_n =y_{n+1}-y_n$.
  \begin{itemize}
  \item[(i)] If \,\, $\mylim f(y_n)\Delta y_n= c$, \, then \, $\mylim \frac 1n
    \int_{y_0}^{y_n}f(u)du = c$.
  \item [(ii)] If \,\, $f(y_k)\Delta y_k\le c$, \,$k\le n$, \,  then \, $\frac 1n
    \int_{y_0}^{y_n}f(u)du \le c$.
  \item [(iii)] If \,\, $f(y_k)\Delta y_k\ge c$, \,$k\le n$, \, then \, $\frac 1n
    \int_{y_0}^{y_n}f(u)du \ge c$.
  \end{itemize}
\end{lemma}

\begin{proof}
  Since in case (i) we also have $\mylim f(y_{n+1})\Delta{y_n} = c$, by
  Toeplitz Lemma we therefore conclude that
  \begin{equation}
    \label{lim:T}
    \frac 1n\sum_{i=0}^{n-1}f(y_{i})\Delta y_i\to c, \quad
    \frac 1n\sum_{i=0}^nf(y_{i+1})\Delta y_i\to c \quad
    \mbox{as} \quad n\to\infty.
  \end{equation}
  If $f$ is decreasing, then by geometrical consideration it is clear that
  \begin{equation}
    \label{eq:2ineq}
    \sum_{i=0}^{n-1}f(y_{i})\Delta
    y_i\ge \int_{y_0}^{y_{n}}f(u)du \ge \sum_{i=1}^nf(y_{i+1})\Delta y_i,
  \end{equation}
  and the result follows from (\ref{lim:T}).  If $f$ is an increasing
  function, then we reverse inequalities in (\ref{eq:2ineq}).

  To prove (ii) in case of decreasing $f$ we note that
  \begin{equation*}
  \int_{y_0}^{y_{n}}f(u)du \le \sum_{i=0}^{n-1}f(y_{i})\Delta y_i \le cn.
  \end{equation*}
  Case of increasing $f$ and (iii) are analogous.
\end{proof}

\begin{corollary}
  \label{cor:gamma}
  Let  $c, \gamma>0$ be non-random numbers.  Let $\{y_n\}_{n\in \mathbb{N}}$
  be a positive increasing sequence and let
\begin{equation}
 \label{cond:ydelt}
  \mylim \Delta y_n =0, \quad  \mylim y_n^{\gamma}\Delta y_n=c.
  \end{equation}
  Then
  \begin{equation*}
  \frac{y_n}{n^{\frac1{1+\gamma}}}\to
  \bigl(c(1+\gamma)\bigr)^{\frac1{1+\gamma}}.
  \end{equation*}
\end{corollary}
\begin{proof}
  We put $f(u)=u^{\gamma}$ for $u\ge 0$, and note that $y_n\to \infty$
  by (\ref{cond:ydelt}). Therefore
  \begin{equation*}
  \mylim \frac{f(y_n)}{f(y_{n-1})}=\mylim
  \left(\frac{y_n}{y_{n-1}}\right)^{\lambda}=1.
  \end{equation*}
  Then by Lemma \ref{lem:T}
  \begin{equation*}
  c=\mylim \frac 1n \int_{y_0}^{y_{n}}f(u)du=\mylim \frac
  {y^{\gamma+1}_{n+1}-y^{\gamma+1}_{0}}{n(\gamma +1)}=\mylim \frac
  {y^{\gamma+1}_{n+1}}{n(\gamma +1)},
  \end{equation*}
  and result follows.
\end{proof}

%Lemma below can be obtained from Toeplitz lemma (Lemma
%\ref{lem:toeplitz}) and the Law of Iterated Logarithm.
%\begin{lemma}
%  \label{lem:g1}
%  Let $\{Y_n\}_{n\in \mathbf{N}}$ be a martingale, $Y_n=\sum_{i=1}^{n}\rho_i$.
%  If a.s. $\mathbf{E}[\rho_n^2\bigr|\mathcal{F}_n]\to 1$ as $n\to \infty$,
%  then a.s.
%  \begin{equation*}
%  \limsup_{n\to \infty}\frac{Y_n}{\sqrt{n}}= \infty, \quad
%  \lim_{n\to \infty}\frac{Y_n}{n}= 0.
%  \end{equation*}
%\end{lemma}

%%%%%%%%%%%%%%%%%%%%%%%%%%%% Subsection %%%%%%%%%%%%%%%%%%%%%%%%%%%%%%%%%%
\subsection{A comparison theorem}

\begin{theorem}
  \label{thm:compar}
  Suppose that $f$ and $g$ are bounded with $f(0)=g(0)=0$ and the random
  variables $\xi_n$ satisfy Assumption~\ref{assum:noise}.  Assume, further,
  that $x_n\to0$ a.s., where $x_n$ is a solution of (\ref{eq0}).
  \begin{enumerate}
  \item[a)] If
    \begin{equation}\label{eq:rate_assm_L}
      \lim_{u\to 0}\frac{f(u)}{g^2(u)}=L, \quad L\in \Reals,
    \end{equation}
    then
    \begin{equation}
      \label{eq:rate_L}
      \mylim\frac{\ln |x_n|}{\sum_{i=1}^n g^2(x_i)} = hL - h/2.
    \end{equation}
  \item[b)] \label{item:rate_assm_inf} If $f(u)\le 0$ in a neighborhood of
    $u=0$ and
    \begin{equation}\label{eq:rate_assm_inf}
      \lim_{u\to 0}\frac{|f(u)|}{g(u)^2}=\infty,
    \end{equation}
    then
    \begin{displaymath}
      \mylim\frac{\ln |x_n|}{\sum_{i=1}^n |f(x_i)|} = -h.
    \end{displaymath}
  \end{enumerate}
\end{theorem}

\begin{remark}
  If the initial value $x_0$ is non-zero and the distributions of $\xi_n$ are
  non-atomic, the solution $x_n$ is a.s.\ non-zero for any $n$.
\end{remark}

\begin{remark}
  Theorem~\ref{thm:nonlinunstab} imposes restrictions on possible values of
  $L$.  To have convergent $x_n$ we must have $L\leq 1/2$.  Then, in
  equation~(\ref{eq:rate_L}), $-h+2Lh$ is non-positive which one would expect
  with $x_n\to0$.
\end{remark}

\begin{proof}
  We apply logarithm to both parts of equation (\ref{eq0}) to obtain the
  following representation of the solution of the recursion
  \begin{equation}\label{eq:lnsum}
    \ln|x_k|=\ln|x_0| + \sum_{i=1}^{k-1}
    \ln\left|1+h f(x_i)+\sqrt{h} g(x_i)\xi_{i+1} \right|.
  \end{equation}
  We set
  \begin{eqnarray*}
    \lambda_{i+1} &=& \ln\left|1+h f(x_i)+\sqrt{h} g(x_i)\xi_{i+1}\right|,\\
    \Phi_i &=& \condex{\lambda_{i+1}},\\
    d_{i+1} &=& \lambda_{i+1} - \Phi_i.
  \end{eqnarray*}
  The expectation $\Phi_i$ can be estimated by Theorem~\ref{thm:my_ito},
  \begin{equation}
    \label{eq:Phi_est}
    \Phi_i = hf(x_i)\big(1+o(1)\big) - h \frac{g^2(x_i)}{2} \big(1+o(1)\big).
  \end{equation}
  On the other hand, it is easy to see that $d_{i+1}$ is
  martingale-difference.  Using Theorem~\ref{thm:my_ito} once again (now with
  $\f = \ln^2|1+x|$) we can estimate
  \begin{eqnarray}
    \nonumber
    \condex{d_{i+1}^2} &=& \condex{\lambda_{i+1}^2} - \Phi_i^2
    = hf(x_i) o(1) + h g^2(x_i) \big(1+o(1)\big) - \Phi_i^2 \\
    \label{eq:char_est}
    &=& h g^2(x_n) \big(1+o(1)\big),
  \end{eqnarray}
  where to get to the final result we used condition~(\ref{eq:rate_assm_L})
  and estimate~(\ref{eq:Phi_est}).  We remind the reader that $o(1)\to0$ as
  $i\to\infty$.

  Part a). Now we would like to apply the Toeplitz Lemma
  (Lemma~\ref{lem:toeplitz}) to the left-hand side of (\ref{eq:rate_L}), where
  $\ln|x_n|$ is expanded as in (\ref{eq:lnsum}).  To apply Toeplitz Lemma we
  need to show that the event
  \begin{displaymath}
    \Lambda=\{\o\in \O: \sum_{i=1}^{\i} g^2(x_i)<\i\}
  \end{displaymath}
  has zero probability.  Since, by assumption, $x_i\to0$ a.s., we have
  $f(x_i)=O(g^2(x_i))$ as $i\to\infty$. Thus we conclude that on (almost all
of)
  $\Lambda$ the series $\sum_{i=1}^{\i} |f(x_i)|$ is convergent too.
  Consequently, the series $\sum_{i=1}^{\i} \Phi_i$, $\sum_{i=1}^{\i}
  \Phi_i^2$ and $\sum_{i=1}^{\i} \condex{d_{i+1}^2}$ are all absolutely
  convergent on $\Lambda$.

  We now rewrite (\ref{eq:lnsum}) as
  \begin{equation*}
    \ln|x_k|=\ln|x_0| + \sum_{i=1}^{k-1} \Phi_i + \sum_{i=1}^{k-1} d_{i+1}
  \end{equation*}
  and notice that the right-hand side is convergent.  This, however,
  contradicts our assumption that $x_i\to0$ a.s.  We conclude that $\Lambda$
  is a zero-probability event.  As a consequence, we obtain that the
  characteristic of the martingale $\sum_{i=1}^n d_{i+1}$ is divergent too,
  see equation~(\ref{eq:char_est}).

  We can now write
  \begin{equation}
    \label{eq:eval_lnx_over_sum_g2}
    \mylim\frac{\ln |x_n|}{\sum_{i=1}^n g^2(x_i)}
    =  \mylim\frac{\sum_{i=1}^{n-1} \Phi_i}{\sum_{i=1}^n g^2(x_i)}
    + \mylim\frac{\sum_{i=1}^{n-1} d_{i+1}}{\sum_{i=1}^n g^2(x_i)}.
  \end{equation}
  The first limit on the right can be evaluated using (\ref{eq:Phi_est}) and
  Lemma~\ref{lem:toepl_cor1},
  \begin{equation*}
    \mylim\frac{\sum_{i=1}^{n-1} \Phi_i}{\sum_{i=1}^n g^2(x_i)}
    = h \lim_{i\to\infty} \frac{f(x_i)\big(1+o(1)\big)}{g^2(x_i)}
    - \frac{h}{2} \lim_{i\to\infty}\frac{g^2(x_i)\big(1+o(1)\big)}{g^2(x_i)}
    = hL - h/2.
  \end{equation*}
  The second limit in~(\ref{eq:eval_lnx_over_sum_g2}) is represented as
  \begin{equation*}
    \mylim\frac{\sum_{i=1}^{n-1} d_{i+1}}{\sum_{i=1}^n g^2(x_i)}
    = \mylim \frac{\sum_{i=1}^{n-1} d_{i+1}}
    {\sum_{i=1}^{n-1}\condex{d_{i+1}^2}}
    \mylim \frac{\sum_{i=1}^{n-1}\condex{d_{i+1}^2}}{\sum_{i=1}^{n} g^2(x_i)}.
  \end{equation*}
  While the second limit on the right is finite (equal to $h$ to be precise),
  the first one is zero by Lemma~\ref{lem:div_variance}.

  Part b).  We follow the proof of part a) with $f(x_i)$ instead of
  $g^2(x_i)$.  By a similar reasoning we conclude that the series
  $\sum_{i=1}^\infty |f(x_i)|$ is divergent.  Indeed, if it were not so, the
  series $\sum_{i=1}^\infty g^2(x_i)$ would be convergent too, since
  $g^2(x_i)/|f(x_i)|\to0$.  We then conclude that the series $\sum_{i=1}^{\i}
  \Phi_i$, $\sum_{i=1}^{\i} \Phi_i^2$ and $\sum_{i=1}^{\i} \condex{d_{i+1}^2}$
  are all absolutely convergent and therefore $\ln|x_k|$ converges to a finite
  limit.  This contradicts our assumptions.

  We write
  \begin{equation}
    \label{eq:eval_lnx_over_sum_f}
    \mylim\frac{\ln |x_n|}{\sum_{i=1}^n |f(x_i)|}
    =  \mylim\frac{\sum_{i=1}^{n-1} \Phi_i}{\sum_{i=1}^n |f(x_i)|}
    + \mylim\frac{\sum_{i=1}^{n-1} d_{i+1}}{\sum_{i=1}^n|f(x_i)|},
  \end{equation}
  and evaluate the first limit using Lemma~\ref{lem:toepl_cor1},
  \begin{equation*}
    \mylim\frac{\sum_{i=1}^{n-1} \Phi_i}{\sum_{i=1}^n g^2(x_i)}
    = h \lim_{i\to\infty} \frac{f(x_i)\big(1+o(1)\big)}{|f(x_i)|}
    - \frac{h}{2} \lim_{i\to\infty}\frac{g^2(x_i)\big(1+o(1)\big)}{|f(x_i)|}
    = -h - \frac{h}{2}\times0.
  \end{equation*}
  To conclude the proof, the second limit in (\ref{eq:eval_lnx_over_sum_f}) is
  evaluated to zero as in part a).
\end{proof}

%%%%%%%%%%%%%%%%%%%%%%%%%%%% Subsection %%%%%%%%%%%%%%%%%%%%%%%%%%%%%%%%%%
\subsection{Rate of decay of $\ln|x_n|$}

Theorem~\ref{thm:compar} provides some information on the decay of
solutions $x_n$ to zero, but does it in a rather implicit way.  We
will now show how one can extract an explicit estimate on the
decay of $\ln|x_n|$ as a function of $n$.

Consider the following example. Let it be given that $x_n\to0$
a.s.\ as $n\to\infty$ (see Theorem~\ref{thm:stability} for a set
of sufficient conditions) and let condition (\ref{eq:rate_assm_L})
be satisfied. Assume that the function $g(u)$ behaves like a power
of $u$ around zero,
\begin{equation*}
  \lim_{u\to0} \frac{g^2(u)}{u^{\mu_g}} = \textrm{const}.
\end{equation*}
Then, using the following lemma we can conclude that
\begin{equation*}
  \mylim \frac{\ln|x_n|}{\ln \left(n^{-1/\mu_g}\right)} = 1.
\end{equation*}

\begin{lemma}\label{lem:ln_invert}
  Let $\lambda>0$ and $x_n$ be a positive sequence satisfying
  \begin{equation}\label{eq:ln_invert_cond}
    \mylim x_n = 0 \quad \mbox{ and } \quad
    \mylim \frac{\ln x_n}{\sum_{i=1}^n x_i^\lambda} = -b < 0.
  \end{equation}
  Then
  \begin{displaymath}
    \mylim \frac{\ln x_n}{\ln n} = -1/\lambda.
  \end{displaymath}
\end{lemma}

\begin{proof}

  We will prove the lemma for $\lambda=1$; the general case would follow with
  the change $x_i \mapsto x_i^\lambda$.  Consider a new variable $y_n =
  \sum_{i=1}^n x_i$.  With $\Delta{y_n} = y_n-y_{n-1}$ we can rewrite the second condition in
  (\ref{eq:ln_invert_cond}) as
  \begin{displaymath}
    \mylim \frac{\ln\Delta{y_n}}{y_n} = -b.
  \end{displaymath}
  By applying the definition of limit,
  given a small $\e>0$ we can find $N$ such that for all $n\geq N$
  \begin{displaymath}
    \frac{\ln\Delta{y_n}}{y_n} \leq -(b-\e) < 0.
  \end{displaymath}
  This can be transformed into
  \begin{displaymath}
    1 \geq e^{(b-\e)y_n}\Delta{y_n}.
  \end{displaymath}
  Applying Lemma \ref{lem:T}, (ii), we obtain
  \begin{equation*}
  \frac1{n(b-\e)} \left[e^{(b-\e)y_n} - e^{(b-\e)y_N}\right]
  = \frac 1n \int_{y_N}^{y_n}e^{(b-\e)u}du\le 1.
  \end{equation*}
  Solving for $y_n$ gives
  \begin{displaymath}
    y_n \leq \frac{\ln\left[C_1n+C_2\right]}{(b-\e)},
  \end{displaymath}
  where $C_1$ and $C_2$ are constant with respect to $n$.  Sending $n$ to
  infinity we obtain
  \begin{equation}\label{eq:lim_bound1}
    \limsup_{n\to\infty} \frac{y_n}{\ln n} \leq \frac{1}{(b-\e)}.
  \end{equation}
  To obtain a similar bound from below we notice that
  \begin{displaymath}
    \frac{y_n}{y_{n-1}} = \frac{y_{n-1}+x_n}{y_{n-1}} \to 1,
  \end{displaymath}
  and therefore
  \begin{displaymath}
    \mylim \frac{\ln\Delta{y_n}}{y_{n-1}} = -b.
  \end{displaymath}
   Retracing our steps with $+\e$ instead of $-\e$  and applying Lemma \ref{lem:T}, (iii),
   produce
\begin{equation*}
 \frac{1}{n(b+\e)} \left[e^{(b+\e)y_n} - e^{(b+\e)y_N}\right]=\frac 1n \int_{y_N}^{y_n}e^{(b+\e)u}du\ge 1
\end{equation*}

  and, ultimately,
  \begin{equation}\label{eq:lim_bound2}
    \liminf_{n\to\infty} \frac{y_n}{\ln n} \geq \frac{1}{(b+\e)}.
  \end{equation}
  Since $\e$ was arbitrary, we obtain from (\ref{eq:lim_bound1}) and
  (\ref{eq:lim_bound2})
  \begin{displaymath}
    \mylim \frac{y_n}{\ln n} = \frac{1}{b},
  \end{displaymath}
  and, using the definition of $y_n$ and condition (\ref{eq:ln_invert_cond}),
  \begin{displaymath}
    \mylim \frac{\ln x_n}{\ln n}
    = \mylim \frac{\ln x_n}{y_n}  \frac{y_n}{\ln n}
    = -b \frac1b = -1.
  \end{displaymath}
\end{proof}

Next we formulate a corollary which extends and formalizes the
discussion at the start of the present section.

\begin{corollary}\label{cor:rate}
  Suppose that $f$ and $g$ are bounded with $f(0)=g(0)=0$ and the random
  variables $\xi_n$ satisfy Assumption~\ref{assum:noise}.  Assume, further,
  that $x_n\to0$ a.s., where $x_n$ is a solution of (\ref{eq0}).
  If one of the following conditions is fulfilled,
  \begin{itemize}
  \item[a)]
    \begin{equation}\label{eq:cor_cond_A}
      \lim_{u\to 0}\frac{f(u)}{|u|^\lambda} = c < 0,\qquad
      \lim_{u\to 0}\frac{f(u)}{ g(u)^2}=-\i,
    \end{equation}
    or
  \item [b)]
    \begin{displaymath}
      \lim_{u\to 0}\frac{g^2(u)}{|u|^\lambda} = c > 0,\qquad
      \lim_{u\to 0}\frac{f(u)}{ g(u)^2} = L < \frac12,
    \end{displaymath}
  \end{itemize}
  then
  \begin{displaymath}
    \mylim\frac{\ln |x_n|}{\ln n}=-\frac{1}{\lambda}.
  \end{displaymath}
\end{corollary}

\begin{proof}
  For the solution $x_n$ we have
  \begin{displaymath}
    \mylim \frac{\ln|x_n|}{\sum_{i=1}^n |x_i|^\lambda}
    = \mylim \frac{\ln|x_n|}{\sum_{i=1}^n f(x_i)}
    \mylim \frac{\sum_{i=1}^n f(x_i)}{\sum_{i=1}^n |x_i|^\lambda}
    = h \mylim \frac{f(x_i)}{|x_i|^\lambda} = hc < 0,
  \end{displaymath}
  where the first limit was calculated using Theorem~\ref{thm:compar} and the
  second was done using Lemma~\ref{lem:toepl_cor1} and
  condition~(\ref{eq:cor_cond_A}).  Now we apply Lemma~\ref{lem:ln_invert} to
  finish the proof of Part a).  The proof of Part b) is analogous.
\end{proof}

\begin{remark}
\label{rem:ln} Relation
\begin{displaymath}
    \mylim \frac{\ln |x_n|}{\ln n} = -\frac 1{\lambda}<0, \quad \mbox{a.s.}
  \end{displaymath}
  implies that for every $\varepsilon>0$ there exists $N=N(\varepsilon,
  \omega)$ such that a.s. for all $n\ge N(\varepsilon, \omega)$
\begin{displaymath}
    n^{-\frac 1{\lambda}-\varepsilon}\le |x_n|\le n^{-\frac 1{\lambda}+\varepsilon}.
  \end{displaymath}
  \end{remark}

We would also like to mention that Lemma~\ref{lem:ln_invert} can
be extended to include other forms of the functions $f(u)$ and
$g^2(u)$ around the origin.  We give this result here without the
proof, which is a simple extension of the proof of
Lemma~\ref{lem:ln_invert} (see also Theorem 5 of \cite{ARS}).

\begin{lemma}
  \label{lem:ln_invert_gen}
  Let $x_n$ be a positive sequence satisfying
  \begin{equation*}
    \mylim x_n = 0 \quad \mbox{ and } \quad
    \mylim \frac{\ln |x_n|}{\sum_{i=1}^n |f(x_i)|} = -b < 0.
  \end{equation*}
  Suppose there exists a function $a(u)$ satisfying
  \begin{itemize}
  \item $a(u)$ is monotone increasing,
  \item $|f(u)|/a(u)\to 1$ as $u\to0+$,
  \item the function $A(z)$, defined by $A(z) = \int_z^1 \frac{du}{ua(u)}$,
    obeys
    \begin{equation*}
      \lim_{y\to\infty}
      \frac{\ln\left[A^{-1}(y-y^*)\right]}{\ln\left[A^{-1}(y)\right]}
      = 1
    \end{equation*}
    for any constant $y^*$.
  \end{itemize}
  Then
  \begin{displaymath}
    \mylim \frac{\ln |x_n|}{\ln\left[A^{-1}(n)\right]} = 1.
  \end{displaymath}
\end{lemma}

%%%%%%%%%%%%%%%%%%%% Polynomial Decay Rate %%%%%%%%%%%%%%%%%%%%%%%%%%%%%%%%%%
\section{Exact Rate of Decay}
\label{sec:exact_rate}

In this section we derive the exact decay rate (or prove its
absence) in the case when $f$ and $g$ have power-law behavior:
\begin{equation}
 \label{eq:polfg}
 f(u)\sim-a_f|u|^{\mu_f}, \quad g^2(u)\sim a_g|u|^{\mu_g} \quad
 \mbox{as} \quad u\to 0,
\end{equation}
where $\mu_f, \mu_g, a_g>0$ and $a_f\neq 0$.

We will assume that
\begin{equation}
  \label{eq:poly_limit}
  \lim_{u\to 0}\frac{2f(u)}{g^2(u)} < 1.
\end{equation}
and that $x_n\to0$.

The assumptions above ensure that the conditions of
Corollary~\ref{cor:rate} are satisfied, which gives us a
preliminary estimate on the rate of decay of $x_n$ (see
Remark~\ref{rem:ln}). The aim of this section is to strengthen the
result of Remark~\ref{rem:ln} to the result of the type
$x_nn^\lambda\to\textrm{const}$ or to show that such strengthening
is impossible.

\begin{remark}
  \label{rem:three_cases}
  It is clear that if $f$ and $g$ are given by (\ref{eq:polfg}) then
  condition~(\ref{eq:poly_limit}) holds in the following three cases:
  \begin{enumerate}
  \item $\mu_f > \mu_g$,
  \item $\mu_f=\mu_g$ and $-2a_f < a_g$,
  \item  $\mu_f < \mu_g$ and $a_f > 0$.
  \end{enumerate}
\end{remark}

We will also need to strengthen Assumption 1 about the noise
$\xi_n$:
\begin{assumption}
  \label{assum:noise1}
  Let $\xi_n$ be independent random variables.
   %with density functions $p_n(\xi)$.
   We assume that $\E\xi_n = 0$, $\E\xi_n^2=1$ and for each $m\in\mathbb{N}$
   there is a $C(m)>0$ such that
\begin{equation}
    \label{eq:moment_decay}
    \mathbf{E}[|\xi_n|^m]\le C(m).
  \end{equation}
\end{assumption}
An important example of the noise satisfying
Assumption~\ref{assum:noise1} are the i.i.d.\ normal $\xi_n$.
%We note that condition~(\ref{eq:density_decay}) of Assumption
%\ref{assum:noise1} ensures that all moments of $|\xi_n|$ are
%bounded uniformly in $n$.
The condition \eqref{eq:moment_decay} implies that the large
fluctuations of $\xi_n$ grow slower than any power.
\begin{lemma}
  \label{lem:xipolyld}
  Suppose that \eqref{eq:moment_decay} holds. Then for every fixed $\varepsilon>0$
  \begin{equation*}
  \lim_{n\to\infty} n^{-\varepsilon}|\xi_n|=0, \quad \text{a.s.}
  \end{equation*}
\end{lemma}

This lemma is a direct consequence of \eqref{eq:moment_decay} and
the Borel-Cantelli lemma.

\begin{corollary}\label{cor:xxi}
  From power-law decay of $x_n$ (Remark~\ref{rem:ln}) and
  Lemma~\ref{lem:xipolyld} we conclude that for every fixed $\varepsilon>0$
  \begin{equation*}
  \lim_{n\to\infty} |x_n|^\varepsilon |\xi_n|=0, \quad \text{a.s.}
  \end{equation*}
  In particular, $g(x_n)\xi_{n+1}\to 0$ as $n\to\infty$.
\end{corollary}

For the sake of simplicity everywhere below we are going to hide
$h$ and $\sqrt{h}$. In other words we let
\begin{equation*}
f:=hf, \quad g:=\sqrt{h}g.
\end{equation*}
\begin{remark}
  \label{rem:pos}
  Since $x_n\to 0$, $f(x_n)\to0$ and $g(x_n)\xi_{n+1}\to 0$ as
  $n\to\infty$, after some random number of steps the bracket
  \begin{equation*}
  \left(1+f(x_n)+g(x_n)\xi_{n+1}\right)
  \end{equation*}
  becomes positive. In particular, it means that solution $x_n$ eventually
  stops changing sign. %Therefore, we may assume without loss of
%  generality that convergence takes place on a nontrivial event on
%  which $x_n>0$ for all $n>N$, where $N$ is an almost surely
%  finite random variable.
\end{remark}

%%%%%%%%%%%%%%%%%%%%%%%%%%%% New Subsection %%%%%%%%%%%%%%%%%%%%%%%%%%%%
\subsection{Main Results}

It turns out that, in the situation described in case (iii) of
Remark~\ref{rem:three_cases}, the decay rate is exact.
\begin{theorem}
  \label{thm:exact}
  Suppose Assumption 2 holds. Let $a_f>0$ and $\mu_f<\mu_g$. Then
  \begin{equation*}
    \lim_{n\to \infty}|x_n|n^{\frac{1}{\mu_f}}=
    \left[\frac{1}{a_f \mu_f}\right]^{\frac 1{\mu_f}},
    \quad\text{a.s.}
  \end{equation*}
\end{theorem}
The intuitive reason for the above behavior is that the
convergence under the conditions of the Theorem is dominated by
the deterministic terms of the recursion and the rate coincides
with the corresponding deterministic ($\xi_n\equiv0$) rate.

On the other hand, if conditions of cases (i) or (ii) of
Remark~\ref{rem:three_cases} are met, then $x_n$ undergoes large
oscillations around the power law decay. This is because the
convergence is induced by the noise term, which is now
significant. More detailed explanations are given in
Remark~\ref{rem:outline_exact} below and, of course, in the proofs
in Sections~\ref{sec:proof_exact} and \ref{sec:proof_nexact}.

\begin{theorem}\label{thm:nexact}
Suppose that Assumption 2 holds. Let either $\mu_f>\mu_g$ or
$\mu_f=\mu_g$ and $-2a_f < a_g$ hold (cases (i)and (ii) of
Remark~\ref{rem:three_cases} correspondingly). Suppose moreover
that $g(u)\ge 0$ for $u>0$ and there exists some $r>0$ such that
\begin{equation}
 \label{cond:fad}
 \frac 1{\sqrt{a_g}}u^{-\mu_g/2}g(u)=1+o(u^{r}) \quad
 \mbox{as} \quad u\to +0.
\end{equation}
%We also suppose in addition that $\xi_n$ are identically
%distributed.
    Then
  \begin{equation}
    \label{eq:nexct1}
    \limsup_{n\to
    \infty}|x_n|n^{\frac{1}{\mu_g}}=\infty,\quad\text{a.s.}
  \end{equation}
  \begin{equation}
    \label{eq:nexct2}
    \liminf_{n\to \infty}|x_n|n^{\frac{1}{\mu_g}}=0, \quad\text{a.s.}
  \end{equation}
\end{theorem}

%%%%%%%%%%%%%%%%%%%%%%% New Subsection %%%%%%%%%%%%%%%%%%%%%%%%%%%%%%
\subsection{Main Construction and an Outline of the Proofs}

By squaring  both part of (\ref{eq0}) with $f:=hf, \quad
g:=\sqrt{h}g$, we obtain the equation
\begin{equation}
 \label{eq:sqr}
  x^2_{n+1}=x^2_n\left(1+F_n+R_{n+1}\right), \quad n=0,1, \dots,
\end{equation}
where
\begin{eqnarray}
 \label{def:Fn}
 F_n&=&2f(x_n)+f^2(x_n)+g^2(x_n),\\
 \label{def:Rn+1}
 R_{n+1}&=&2(1+f(x_n))g(x_n)\xi_{n+1}+g^2(x_n)(\xi^2_{n+1}-1).
\end{eqnarray}
\begin{remark}
 \label{rem:FR0}
Under conditions of this section, $ F_n, R_{n+1} \to 0$ a.s. as
$n\to \infty$.
\end{remark}

 Let $a,\mu\in \mathbb{R}$, $\mu>0$ and $\nu=\frac{\mu}{2}$. We
 define for $u>0$
\begin{equation}
  \label{def:G}
  G(u) = \frac{u^{-\nu}}{a\nu}, \quad x>0.
\end{equation}
Assuming that $x_n\neq 0$, we apply the Taylor expansion with
three terms to obtain for $x^2_n$ a.s.
\begin{equation}
  \label{eq:G}
  \begin{split}
  G(x^2_{n+1}) &= G\Bigl(x^2_n+x^2_n(F_n+R_{n+1})\Bigr)\\
  &= G(x^2_n)-\frac 1a |x_n|^{-2\nu}(F_n+R_{n+1})
  + \frac{\nu+1}{2a}\eta_{n+1}^{-\nu-2}|x_n|^{4}(F_n+R_{n+1})^2,
  \end{split}
\end{equation}
with
\begin{equation}
  \label{eq:kappa}
  |\eta_{n+1}-x^2_n|\le x^2_n|F_n+R_{n+1}|.
\end{equation}
Let
\begin{equation*}
  \kappa_{n+1}=\eta_{n+1}^{-(\nu+2)}x_n^{2(\nu+2)}.
\end{equation*}

We can rewrite (\ref{eq:G}) in the following form
\begin{multline}
  \label{eq:G1}
  %\begin{split}
  G(x^2_{n+1})= G(x^2_n)-\frac 1a |x_n|^{-\mu}(F_n+R_{n+1})
  +\frac{\nu+1}{2a}|x_n|^{-\mu}(F_n+R_{n+1})^2\\
  +\frac{\nu+1}{2a}(\kappa_{n+1}-1)|x_n|^{-\mu}(F_n+R_{n+1})^2\\
 =G(x^2_n) + P_{n+1} + \rho_{n+1} + \tau_{n+1},
 %\end{split}
\end{multline}
where
\begin{align}
  \label{def:P}
  P_{n+1} &= -\frac2a |x_n|^{-\mu}\left[f(x_n)
    - \frac{\mu+1}{2}g^2(x_{n})\right] + Q_{n+1},\\
  \label{def:rho1}
  \rho_{n+1}&=-\frac 2a |x_n|^{-\mu}g(x_{n})\xi_{n+1},\\
  \label{def:rho2}
  \tau_{n+1}&=\frac {\mu+1}{a}|x_n|^{-\mu}g^2(x_{n})(\xi^2_{n+1}-1)
  ,\\
  \label{def:Q}
  Q_{n+1} =&-\frac{1}{a}|x_n|^{-\mu}\Big(f^2(x_n)+2f(x_n)g(x_{n})\xi_{n+1}\Big)\nonumber\\
&+\frac{\mu+2}{4a}|x_n|^{-\mu}\Big(F_n^2+2F_nR_{n+1}+(\kappa_{n+1}-1)(F_n+R_{n+1})^2\Big)\\
&+\frac{\mu+2}{4a}|x_n|^{-\mu}\Big(4f^2(x_n)g^2(x_{n})\xi^2_{n+1}+8f(x_n)g^2(x_{n})\xi^2_{n+1}\nonumber\\
&+g^4(x_{n})(\xi^2_{n+1}-1)^2+4g^3(x_{n})(1+f(x_n))(\xi^2_{n+1}-1)\xi_{n+1}\Big)
.\nonumber
\end{align}
Let
\begin{equation}
  \label{def:mu}
  \mu=\left\{
    \begin{array}{cc}
      \mu_f &\mbox{if} \quad  \mu_f<\frac{\mu_g}2,\\
      \mu_g-\mu_f & \qquad \mbox{if} \quad  \frac{\mu_g}2\le \mu_f< \mu_g, \\
      \frac{\mu_g}2 &\mbox{if} \quad  \mu_f> \mu_g. \\
  \end{array}
\right.
\end{equation}
Since $\rho_n$ and $\tau_n$ defined by
(\ref{def:rho1})-(\ref{def:rho2}) are
$\mathcal{F}_n$-martingale-differences,
\begin{equation}
  \label{def:Mi}
  M_n=\sum_{i=0}^{n-1}\rho_{i+1} \qquad\mbox{and}\qquad
  T_n=\sum_{i=0}^{n-1}\tau_{i+1}
\end{equation}
are $\mathcal{F}_n$-martingales. After summation of (\ref{eq:G1}),
we arrive at
\begin{equation}
  \label{eq:Gsum}
  G(x^2_{n})=G(x^2_0)+\sum_{i=0}^{n-1}P_{i+1}+M_n+T_n.
\end{equation}

\begin{remark}{\bf (Outline of the proof of Theorems~\ref{thm:exact} and
    \ref{thm:nexact}).}
\label{rem:outline_exact}
  The decomposition of $G(x^2_{n+1})$, equation (\ref{eq:G1}), is constructed in
  a way that highlights the different types of convergence as $n \to \infty$
  for different values of $\mu$.

  Firstly, the term $P_{n+1}$ is made up from 3 parts.  The last part,
  $Q_{n+1}$, is subdominant to the first two for all values of $\mu$.  When
  $\mu_f<\mu_g$ the limiting behavior of $P_{n+1}$ is dominated by the
  function $f$, while when $\mu_f>\mu_g$, the function $g$ determines the
  behavior.

  When $\mu_f<\mu_g/2$ we set $\mu=\mu_f$ and show that both
  $\rho_{n+1}$ and $\tau_{n+1}$ a.s.\ tend to zero as $n\to
  \infty$ and $P_{n+1}$ tends to 1 as $n\to \infty$.  This means that
  asymptotic behavior of $G(x^2_{n+1})-G(x^2_{n})$ is determined by $P_{n+1}$ and
  the result can be obtained directly by Toeplitz lemma.

  When $\mu_g>\mu_f\ge \mu_g/2$, the martingale-difference $\rho_{n+1}$ no
  longer tends to zero if $\mu=\mu_f$.  To get around this difficulty we set
  $\mu=\mu_g-\mu_f$ and show that in this situation both $P_{n+1}+\tau_{n+1}$
  and $\condexn{\rho_{n+1}^2}$ behave like $|x_n|^{2\mu_f-\mu_g}$ as $n\to
  \infty$.  By means of a martingale convergence theorem, namely by Lemma
  \ref{lem:div_variance}, we compare the behavior of $G(x^2_{n})$ with that of
  $\sum_{i=1}^{i=n}|x_i|^{2\mu_f-\mu_g}$ and apply Corollary \ref{cor:gamma}.

  When $\mu_f\ge \mu_g$, both $\rho_{n+1}$ and $\tau_{n+1}$ decay slower
  than $P_{n+1}$ for all values of $\mu$.  We set $\mu=\mu_g/2$ so that
  $\condexn{\left(\rho_{n+1}+\tau_{n+1}\right)^2}\to1$ and $P_{n+1}\sim
  |x_n|^{\mu_g/2}$ as $n\to\infty$.  This allows us to apply a consequence of the central limit theorem
  (see Lemma~\ref{lemma:clt_liapunov} below) and decomposition (\ref{eq:Gsum}) to
  obtain the conclusion of Theorem~\ref{thm:nexact}.
\end{remark}

The rest of this section we devote to verifying that $Q_{n+1}$ is
subdominant to the other terms in $P_{n+1}$.

\begin{lemma}\label{lem:kappa}
  There is some $K=K(\mu)>0$ and $N=N(K, \omega)$ such that for all
  $n\ge N$ we have
  \begin{equation}
    \label{est:kappa}
    |\kappa_{n+1}-1|\le K|F_n+R_{n+1}|.
  \end{equation}
  In particular, a.s. $\kappa_n\to 1$ as $n\to \infty$.
\end{lemma}

\begin{proof}
  From (\ref{eq:kappa}) we have
  \begin{equation}
    \label{est:kappa1}
    \left|\frac{\eta_{n+1}}{x^2_n}-1\right|\le |F_n+R_{n+1}|.
  \end{equation}
  Since for any $\lambda\in \mathbb{R}$,
  \begin{equation*}
  \lim_{y\to 1}\frac{y^{\lambda}-1}{y-1}=\lambda,
  \end{equation*}
  for any $K>|\lambda|$ we can find $\delta>0$ such that
  \begin{equation}
    \label{est:kappa2}
    |y^{\lambda}-1|\le K|y-1|
  \end{equation}
  when $|y-1|<\delta$. Letting $\lambda=-(\mu/2+2)$, and using
  (\ref{est:kappa1}) and Remark \ref{rem:FR0}, we can find $N(K, \omega)$ such that for all
  $n>N(K, \omega)$
  \begin{equation*}
  |F_n+R_{n+1}|<\delta,
  \end{equation*}
  which, together with (\ref{est:kappa2}), implies
  (\ref{est:kappa}).
\end{proof}

\begin{lemma}
  \label{lem:Q}
  Let $\mu$ be as defined in (\ref{def:mu}), and $s=\min\{\mu_f, \mu_g\}$.
  Then, for any $\varepsilon>0$, a.s. $Q_{n+1}\sim O(x_n^{s-\varepsilon})$ as
  $n\to \infty$.
\end{lemma}

\begin{proof}
  By application of Corollary \ref{cor:xxi} for any $\varepsilon>0$
  we obtain that a.s.
  \begin{equation*}
    \label{est:k2}
    \kappa_{n+1}-1=O\left(|x_n|^{\mu_f}+|x_n|^{\frac{\mu_g}2-\varepsilon}\right),
    \quad \mbox{as} \quad n\to \infty,
  \end{equation*}
  and then, from (\ref{def:Q}),
  \begin{equation}
    \label{est:O2}
    Q_{n+1}=O\left(|x_n|^{-\mu+2\mu_f} + |x_n|^{-\mu+\mu_f+\frac{\mu_g}2-\varepsilon}
     + |x_n|^{-\mu+\frac{3\mu_g}2-\varepsilon}\right).
  \end{equation}
  Now the proof can be completed  by direct substitution of
  different values of $\mu$ from (\ref{def:mu}) into (\ref{est:O2}).
\end{proof}

%%%%%%%%%%%%%%%%%%%%%%%%%%% New Subsection %%%%%%%%%%%%%%%%%%%%%%%%%%%%%%
\subsection{Proof of Theorem \ref{thm:exact}}
\label{sec:proof_exact}

%%%%%%%%%%%%%%%%%%%%%%%%%%% new subsubsection %%%%%%%%%%%%%%%%%%%%%%%%%%%
\subsubsection{An auxiliary lemma}

\begin{lemma}
  \label{lem:Pmuf}
  Let $a_f>0$ and $\mu_f<\mu_g$.  If we set $\mu=\mu_f$, and $a=2a_f>0$ in
  equations (\ref{def:P}), (\ref{def:rho2}) then
  \begin{equation*}
    \lim_{n\to \infty}\bigl[P_{n+1}+\tau_{n+1}\bigr]=1, \quad \text{a.s.}.
  \end{equation*}
\end{lemma}
The result follows from equations (\ref{eq:polfg}), (\ref{def:P}),
(\ref{def:rho2}), Corollary \ref{cor:xxi} and Lemma \ref{lem:Q}.

%%%%%%%%%%%%%%%%%%%%%%%%%%% new subsubsection %%%%%%%%%%%%%%%%%%%%%%%%%%%
\subsubsection{Proof of Theorem~\ref{thm:exact}}

We consider two cases:
\begin{enumerate}
  \item[(i)] $\mu_f<\frac 12 \mu_g$,
  \item[(ii)] $\frac 12 \mu_g\le \mu_f<\mu_g$.
\end{enumerate}

%%%%%%%%%%%%%%%%%%% Case (i)
\begin{proof}[Proof of Theorem~\ref{thm:exact}, case (i)]
  We set $\mu=\mu_f$, and $a=2a_f>0$ in equations (\ref{def:P})-(\ref{def:Q}).  By Corollary \ref{cor:xxi} for
  $\varepsilon< \frac{\mu_g}2-\mu_f$ we have a.s.  as $n\to \infty$
  \begin{equation*}
  \rho_{n+1}=O\left(|x_n|^{-\mu_f+\frac{\mu_g}2-\varepsilon}\right)\to 0.
  \end{equation*}
  Then  from Lemma \ref{lem:Pmuf} we obtain that for $n\to \infty$
  \begin{equation*}
  G(x^2_{n+1})-G(x^2_n)\to 1, \quad \text{a.s.}
  \end{equation*}
  which together with Toeplitz Lemma and (\ref{def:G}), implies that
  a.s. for $n\to \infty$
  \begin{equation*}
  \frac 1{na_f\mu_f}|x_{n}|^{-\mu_f} = \frac{G(x^2_{n})}{n}\to
  1.
  \end{equation*}
\end{proof}

%%%%%%%%%%%%%%%%%%%% Case (ii)

\begin{proof}[Proof of Theorem~\ref{thm:exact}, case (ii)]
  We set $\mu=\mu_g-\mu_f$, $a=\frac{2a_g}{a_f}$ in equations (\ref{def:P})-(\ref{def:Q}).  We also denote
  \begin{equation*}
  b:=\frac{a^2_f}{a_g}>0, \quad \lambda:=2\mu_f-\mu_g\ge 0.
  \end{equation*}
  Applying again Corollary \ref{cor:xxi} and Lemma \ref{lem:Q} and reasoning
  in the usual way we get a.s.\ as $n\to\infty$
  \begin{equation}
    \label{eq:ratdM}
    \frac{\mathbf{E}[\rho_{n+1}^2|\mathcal{F}_n]}{b|x_n|^{\lambda}}=
    \frac{g^2(x_n)|x_n|^{-2(\mu_g-\mu_f)}}{a_g |x_n|^{2\mu_f-\mu_g}} \to 1,
  \end{equation}
  \begin{equation}
    \label{eq:ratPn}
    \frac{P_{n+1}+\tau_{n+1}}{b|x_n|^{\lambda}}\to 1.
  \end{equation}

  Now we prove that a.s.
  \begin{equation}
    \label{eq:infa}
    \lim_{n\to\infty}\sum_{i=1}^{n}|x_i|^{\lambda}=\infty.
  \end{equation}
  Indeed, let $\Omega_1=\{\omega:\ \lim_{n\to
    \infty}\sum_{i=1}^{n}|x_i|^{\lambda}<\infty\}$. Then from
  (\ref{eq:ratdM}) and (\ref{eq:ratPn}) we obtain that on $\Omega_1$
  we also have
  \begin{equation*}
  \sum_{i=1}^{\infty}[P_{i+1}+\tau_{i+1}]<\infty \qquad\mbox{and}\qquad
  \sum_{i=1}^{\infty} \mathbf{E}[\rho_{i+1}^2|\mathcal{F}_i] < \infty.
  \end{equation*}
  The last relation implies that $\lim_{n\to \infty}M_n$ is a.s.\ finite
  on $\Omega_1$ and then equation (\ref{eq:Gsum}) implies that $\lim_{n\to
    \infty}G(x_n)<\infty$ a.s.\ on $\Omega_1$.  But since $x_n\to0$ a.s., the
  probability of $\Omega_1$ must be zero.

  Relation (\ref{eq:ratdM}) together with (\ref{eq:infa}) implies that a.s.
  \begin{equation*}
    \frac{\sum_{i=1}^n \mathbf{E}[\rho_{i+1}^2|\mathcal{F}_i]}
    {\sum_{i=1}^n b|x_i|^\lambda} \to 1,
  \end{equation*}
  which results in
  \begin{equation}
    \label{eq:ratdM2}
    \frac{M_n}{\sum_{i=1}^nb |x_i|^{\lambda}}\to 0.
  \end{equation}
  On the other hand, relation (\ref{eq:ratPn}) together with (\ref{eq:infa})
  implies that
  \begin{equation*}
    \frac{\sum_{i=1}^n\left[P_{i+1}+\tau_{i+1}\right]}
    {\sum_{i=1}^nb |x_i|^\lambda}\to 1,
  \end{equation*}
  which together with (\ref{eq:ratdM2}) gives
  \begin{equation*}
    \lim_{n\to\infty}\frac{G(x^2_n)}{\sum_{i=1}^n b|x_i|^\lambda}=1.
  \end{equation*}
  After applying (\ref{def:G}) and rearranging, we arrive at
  \begin{equation}
    \label{eq:ratG1}
    |x_n|^{\mu}\sum_{i=1}^n|x_i|^{\lambda}\to \frac 1{a_f\mu}.
  \end{equation}
  Now we are going to apply Corollary \ref{cor:gamma} from Lemma \ref{lem:T}.
  We set $\gamma=\lambda/\mu$ and $c=\left(\frac 1{a_f\mu}\right)^\gamma$.
  We define $y_n:=\sum_{i=1}^{n-1}|x_i|^{\lambda}$, so that $\Delta
  y_n=|x_n|^\lambda$ and $|x_n|^\mu = \left(\Delta y_n\right)^{1/\gamma}$.  Now
  relation (\ref{eq:ratG1}) takes the form
  \begin{equation*}
    \left(\Delta y_n\right)^{1/\gamma}y_{n+1} \to \frac1{a_f\mu},
  \end{equation*}
  which, together with $\frac{y_{n+1}}{y_n}\to 1$, implies
\begin{equation*}
    \left(\Delta y_n\right)^{1/\gamma}y_{n} \to \frac1{a_f\mu},
  \end{equation*}
  or, equivalently,
  \begin{equation}
    \label{eq:ratG3}
    y_n^{\gamma} \Delta y_n \to c.
  \end{equation}
  By applying Corollary \ref{cor:gamma} we obtain
  \begin{equation*}
  \frac{y_n}{n^{\frac1{1+\gamma}}}
  \to \bigl(c(1+\gamma)\bigr)^{\frac1{1+\gamma}},
  \quad\mbox{or}\quad
  \frac{y_n^{\gamma}}{n^{\frac{\gamma}{1+\gamma}}}\to
  \bigl(c(1+\gamma)\bigr)^{\frac{\gamma}{1+\gamma}}.
  \end{equation*}
  The last limit  together with (\ref{eq:ratG3}) gives
  \begin{equation}
    \label{calc:ylim}
    \Delta y_nn^{\frac{\gamma}{1+\gamma}}=\left(y_n^{\gamma}\Delta y_n\right)
    \times \frac{n^{\frac{\gamma}{1+\gamma}}}{y_n^{\gamma}}\to
    c\bigl(c(1+\gamma)\bigr)^{-\frac{\gamma}{1+\gamma}}.
  \end{equation}
  Substituting the values for $ \Delta y_n$, $\gamma$ and $c$ in
  (\ref{calc:ylim}) gives the desired result.
\end{proof}

%%%%%%%%%%%%%%%%%%%%%%%%%%%% New Subsection %%%%%%%%%%%%%%%%%%%%%%%%%%%%%
\subsection{Proof of Theorem \ref{thm:nexact}}
\label{sec:proof_nexact}

%%%%%%%%%%%%%%%%%%%%%%%%%%%% new subsubsection %%%%%%%%%%%%%%%%%%%%%%%%%%
\subsubsection{Auxiliary lemmas}

%\begin{lemma}(see  Durrett and also Shiryaev)
%If $f:(1,\infty)\to(0,\infty)$ is an increasing function such that
%$\int_{1}^\infty \frac{1}{f^2(t)}\,dt <\infty$, $M_n$ is a
%martingale and $\langle M_n\rangle\to \infty$, then
%\[
%\lim_{n\to\infty} \frac{M_n}{f(\langle M\rangle_n)}=0, \quad
%\text{a.s.}
%\]
%\end{lemma}
The following lemma can be considered as a corollary of a version
of strong law of large numbers for square-integrable martingales
(see e.g. \cite{Shiryaev96}, page 519).
\begin{lemma}
\label{lemma:durrettshiryaev} If $M_n$ is a square-integrable
martingale with the quadratic characteristic $\langle M_n\rangle$
and $\langle M_n\rangle\to \infty$, then for any $\gamma>0$
 \[
\lim_{n\to\infty} \frac{M_n}{(\langle M_n\rangle)^{1/2+\gamma}}=0,
\quad \text{a.s.}
\]
\end{lemma}
Define $\Phi\in C(\mathbb{R};\mathbb{R})$ by
\begin{equation} \label{eq:normaldf}
 \Phi(x):=\frac{1}{\sqrt{2\pi}}\int_{-\infty}^xe^{-\frac{1}{2}y^2}\,dy.
\end{equation}
Assumption \ref{assum:noise1}, due to Liapunov, gives rise to the
following form of the Central Limit Theorem (see e.g. \cite{Bil},
page 362).
\begin{lemma}\label{lemma:clt_liapunov}
Let Assumption 2 hold. Then
\[
\lim_{n\to\infty} \mathbb{P}\left[\frac{1}{\sqrt{n}}\sum_{i=1}^n
\xi_i >x\right]=1-\Phi(x), \quad\text{for all $x\in\mathbb{R}$},
\]
where $\Phi$ is given by \eqref{eq:normaldf}.
\end{lemma}
The following result is then a simple adaptation of the argument
presented on p.380--1 in~\cite{Shiryaev96}.
\begin{lemma}\label{lemma:mngeqsqrtn}
Suppose that $\{\xi_n\}_{n\in \mathbb{N}}$ obeys Assumption 2.
Then
\begin{equation} \label{eq:clt2} \limsup_{n\to\infty}
\frac{1}{\sqrt{n}}\sum_{i=1}^n \xi_i = \infty, \quad
\liminf_{n\to\infty} \frac{1}{\sqrt{n}}\sum_{i=1}^n \xi_i =
-\infty, \quad\text{a.s.} \end{equation}
\end{lemma}
\begin{proof}
For $c>0$ define the events
\[
A_c=\{\omega: \limsup_{n\to\infty} \frac{1}{\sqrt{n}} \sum_{i=1}^n
\xi_i>c\},\quad A=\{\omega: \limsup_{n\to\infty}
\frac{1}{\sqrt{n}}\sum_{i=1}^n \xi_i=\infty\}.
\]
Then $A_c\downarrow A$ as $c\to\infty$. The events $A_c$ are tail
events; therefore, by independence of the sequence
$\{\xi_n\}_{n\in \mathbb{N}}$ and the Zero-One Law, it follows
that
\begin{equation} \label{eq:eqAc}
\mathbb{P}[A_c]>0 \text{ for every $c>0$}
\end{equation}
implies $\mathbb{P}[A_c]=1$, and so
$\mathbb{P}[A]=\lim_{c\to\infty} \mathbb{P}[A_c]=1$. Therefore it
suffices to prove \eqref{eq:eqAc} to establish the first part of
\eqref{eq:clt2}.

Using (i) the fact that for any sequence of random variables
$\{\eta_n\}_{n\in \mathbb{N}}$ we have
\[
\{\omega:\limsup_{n\to\infty} \eta_n(\omega)>x\}\supseteq
\{\omega:\eta_n(\omega)>x\text{ i.o.}\}, \quad \text{for all
$x\in\mathbb{R}$},
\]
(ii) the fact that $\mathbb{P}[B_n \text{ i.o.}]\geq
\limsup_{n\to\infty} \mathbb{P}[B_n]$ for any sequence of events
$\{B_n\}_{n\in \mathbb{N}}$, and then
Lemma~\ref{lemma:clt_liapunov} in turn, we get
\begin{align*}
\mathbb{P}[A_c]&=\mathbb{P}\left[\limsup_{n\to\infty}
\frac{1}{\sqrt{n}}\sum_{i=1}^n \xi_i>c\right] \geq
\mathbb{P}\left[ \frac{1}{\sqrt{n}}\sum_{i=1}^n \xi_i>c \text{
i.o.}\right] \\&\geq \limsup_{n\to\infty} \mathbb{P}\left[
\frac{1}{\sqrt{n}}\sum_{i=1}^n \xi_i>c\right] =\lim_{n\to\infty}
\mathbb{P}\left[ \frac{1}{\sqrt{n}}\sum_{i=1}^n \xi_i>c\right]
=1-\Phi(c),
\end{align*}
proving \eqref{eq:eqAc}.

The second part follows from the first using the change
$\xi_i\rightarrow -\xi_i$.
\end{proof}
\begin{lemma}
  \label{lem:M}
  Let  $a=2\sqrt{a_g}$ and $\mu=\mu_g/2$. Then
  \begin{equation}
 \label{lim:1}
 \limsup_{t\to \infty}\frac{-(M_n+T_n)}{\sqrt{n}}=\infty, \quad \limsup_{t\to
 \infty}\frac{M_n+T_n}{\sqrt{n}}=\infty,\quad \mbox{a.s.}
\end{equation}
%and
%  \begin{equation*}
%    \lim_{n\to \infty}\frac{ M_n+T_n}{n}=0 \quad \mbox{a.s.}
%      \end{equation*}
\end{lemma}
%For the the proof it is enough to verify that a.s.
%$\condexn{(\rho_n+\tau_n)^2}\to1$ as $n\to\infty$ and apply Lemma
%\ref{lem:g1}.

\begin{proof}
We consider the case when $\mu_g<\mu_f$, and therefore we are
under conditions of Corollary 2, b). Then, according to Remark 5,
for any $\varepsilon>0$
\begin{equation}
\label{*}
 n^{-\frac 1{\mu_g}-\varepsilon}\le |x_n|\le  n^{-\frac 1{\mu_g}+\varepsilon}
\end{equation}
for all $n\ge N(\varepsilon, \omega)$. We choose
$\varepsilon<\frac 1{\mu_g}$ and in the following consider only
$n\ge N(\varepsilon, \omega)$.

We define $M_n$ and $\rho_{n+1}$ as above and rearrange
\[
\rho_{n+1}:=-\xi_{n+1}+\bar\rho_{n+1},
\]
where
\[
 \bar \rho_{n+1}= \left[1-\frac
 1{\sqrt{a_g}}|x_n|^{-\mu_g/2}g(x_{n})\right]\xi_{n+1}=o(|x_{n}|^{r})\xi_{n+1}.
\]
 Let $\bar M_n=\sum_{i=1}^n \bar \rho_{i}$. We prove that
\[
 \lim_{n\to\infty}\frac{\bar M_n}{\sqrt{n}}=0, \quad
 \lim_{n\to\infty}\frac{T_n}{\sqrt{n}}=0,
\]
which in conjunction with the conclusion of
Lemma~\ref{lemma:mngeqsqrtn} ensures the desired result.
%\[
%\mathbb{P}\left\{\limsup_{n\to\infty}
%\frac{\sum_{i=1}^n\xi_i}{\sqrt{n}}=\infty\right\}=1,
%\]
%then desired result will follow.

We want to estimate $\langle\bar \rho_n\rangle$. Since
\[
 |x_n|^r\le n^{r\left(-\frac1{\mu_g}+\varepsilon\right)}
\]
we have
\[
\langle\bar \rho_n\rangle\le |x_{n}|^{2r}\le K
n^{2r\left(-\frac1{\mu_g}+\varepsilon\right)}<n^{-\delta}
\]
for some $\delta\in (0, 1)$. Since, as $n\to \infty$,
\[
 \sum_{i=1}^n i^{-\delta}\sim \frac 1{1-\delta}n^{1-\delta},
\]
we have
\[
 \langle\bar M_n\rangle\le K\sum_{i=1}^n i^{-\delta}\sim
 K_1n^{1-\delta}.
\]
For fixed $\delta$ we choose $\gamma>0$ such that $
(1/2+\gamma)(1-\delta)<1/2$. Indeed, it is possible for $\gamma
<\frac{\delta}{2(1-\delta)}$. Then
\[
\lim_{n\to\infty} \frac{\left(\langle\bar M_n\rangle\right)
 ^{1/2+\gamma}}{\sqrt{n}}=0, \quad\text{a.s.}
\]
Hence by applying Lemma~\ref{lemma:durrettshiryaev} we have
\begin{equation*}
\begin{split}
 \lim_{n\to \infty}\frac{\bar M_n}{\sqrt{n}}
 &=\lim_{n\to \infty}\frac{\bar M_n}{\left(\langle\bar M_n\rangle\right)^{1/2+\gamma}}
 \lim_{n\to \infty}\frac{\left(\langle\bar M_n\rangle\right)
 ^{1/2+\gamma}}{\sqrt{n}}=0.
%\times\lim_{n\to \infty}\frac{\left(\langle\bar M_n\rangle\right)
% ^{1/2+\gamma}}{\sqrt{n}}\\
% &\le
% 0\times K_2\lim_{n\to \infty}\frac{n^{(1/2+\gamma)(1-\delta)}}{\sqrt{n}}=0.
 \end{split}
\end{equation*}
%In the same way we prove that
%\[
% \lim_{n\to \infty}\frac{-\bar M_n}{\sqrt{n}}=0.
%\]
We define $T_n$ and $\tau_{n+1}$ as before. Since
\[
 \tau_{n+1}=-\frac {\mu_g+2}{4}\sqrt{a_g}
 |x_n|^{\mu_g/2}\left[1+o(|x_n|^{r})\right]^2(\xi^2_{n+1}-1),
\]
applying (\ref{*}) with $\varepsilon<\frac1{2\mu_g}$ we estimate
for $n\ge N(\varepsilon, \omega)$
\begin{equation*}
\begin{split}
 \langle\tau_{n+1}\rangle&=\left(\frac {\mu_g+2}{4}\right)^2a_g|x_n|^{\mu_g}\left[1+o(|x_n|^{r})\right]^4
 (\mathbf{E}\xi^4_{n+1}-1)\\
 &\le K_2|x_n|^{\mu_g}(\mathbf{E}\xi^4_{n+1}-1)
 \le K_3(\omega)n^{-1+\varepsilon\mu_g}(\mathbf{E}\xi^4_{n+1}-1).
\end{split}
\end{equation*}
Then, as $n\to \infty$, because the fourth moments of $\xi_n$ are
uniformly bounded in $n$,
\[
 \langle T_{n}\rangle\le \sum_{i=1}^n
 K_3(\omega)i^{-1+\varepsilon\mu_g}(\mathbf{E}\xi^4_{i+1}-1)\leq
 \sup_{i\in\mathbb{N}}\{\mathbf{E}\xi_i^4-1\}\times K_3(\omega)\sum_{i=1}^ni^{-1+\varepsilon\mu_g}\sim K_4
 n^{\varepsilon\mu_g}.
\]
Therefore, as $2\mu_g \varepsilon<1$, we have $\langle
T_n\rangle/\sqrt{n}\to0$ as $n\to\infty$ a.s. Also, by the strong
law of large numbers for martingales $T_n/\langle T_n\rangle\to0$
as $n\to\infty$, a.s., so
\begin{equation*}
\begin{split}
 \lim_{n\to \infty}\frac{T_n}{\sqrt{n}}&=\lim_{n\to \infty}\frac{T_n}{\langle T_n\rangle}
 \lim_{n\to \infty}\frac{\langle T_n\rangle}{\sqrt{n}}=0.
% 0\times\lim_{n\to \infty}\frac{\langle\bar T_n\rangle}{\sqrt{n}}\\&\le
% 0\times K_4\lim_{n\to \infty}\frac{n^{\varepsilon\mu_g}}{\sqrt{n}}=0.
 \end{split}
\end{equation*}

Applying Lemma~\ref{lemma:mngeqsqrtn} %we have
%\[
%\mathbb{P}\left\{\limsup_{n\to\infty}
%\frac{\sum_{i=1}^n\xi_i}{\sqrt{n}}=\infty\right\}=1,
%\]
%so
a.s. we get
\begin{equation*}
\begin{split}
 \limsup_{n\to\infty}\frac{M_n+T_n}{\sqrt{n}}&=\limsup_{n\to\infty}\frac{\sum_{i=1}^n(-\xi_i)}{\sqrt{n}}
 +\lim_{n\to\infty}\frac{\bar
 M_n}{\sqrt{n}}+\lim_{n\to\infty}\frac{T_n}{\sqrt{n}}=\infty,
 \end{split}
\end{equation*}
as required. The proof of the 2nd part of (\ref{lim:1}) is
similar.
\end{proof}

\begin{lemma}
  \label{lem:Pmug}
  Let either $\mu_f>\mu_g$ or $\mu_f=\mu_g$ and $-2a_f < a_g$ hold.  We set
  $\mu=\mu_g/2$, and $a=2\sqrt{a_g}>0$ in equation~(\ref{def:P}).  Then
  \begin{equation}
    \label{eq:Pn}
    \lim_{n\to \infty}\frac{P_{n+1}}{a(\mu+1)|x_n|^{\mu}}=S>0,
  \end{equation}
  where the number $S$ is non-random.
\end{lemma}
\begin{proof}
  From equation~(\ref{eq:polfg}) and Lemma \ref{lem:Q} we
  obtain that a.s. as $n\to \infty$
  \begin{equation*}
    \begin{split}
      -\frac 2{a^2(\mu+1) |x_n|^{2\mu}}f(x_n)&\to L=\left\{ \begin{array}{cc}
          0, &  \mu_f>\mu_g,\\
          \frac{a_f}{a_g(\mu_g+2)}, & \mu_f=\mu_g, \\
        \end{array}\right.\\
      \frac {g^2(x_n)}{a^2|x_n|^{2\mu}}&\to \frac 14, \\
      \frac {Q_{n+1}}{a(\mu+1) |x_n|^{\mu}}&\to 0.
    \end{split}
  \end{equation*}
  The above relations imply (\ref{eq:Pn}) with $S=L+\frac 14$.  Note that
  $S>0$ even if $a_f$ is negative.  Indeed, if this is the case and
  $\mu_f=\mu_g$
  \begin{equation*}
  L = \frac{a_f}{a_g(\mu_g+2)} > \frac{-a_g}{2a_g(\mu_g+2)}
  = -\frac{1}{2(\mu_g+2)} > -\frac 14,
  \end{equation*}
  since $-2a_f < a_g$.  The lemma is proved.
\end{proof}

%%%%%%%%%%%%%%%%%%%%%%%% New subsubsection %%%%%%%%%%%%%%%%%%%%%%%%%%
\subsubsection{Proof of Theorem \ref{thm:nexact}.}

We set $\mu=\mu_g/2$, $a=2\sqrt{a_g}$.

\begin{proof}[Proof of (\ref{eq:nexct1})]
  First we rearrange (\ref{eq:Gsum}) in the following way
  \begin{equation*}
    \sum_{i=0}^{n-1}P_{i+1}=G(x_{n})-G(x_0)-M_n-T_n.
  \end{equation*}
  Due to Lemma \ref{lem:M},
  \begin{equation*}
  \limsup_{n\to \infty}\frac{-(M_n+T_n)}{\sqrt{n}}= \infty,
  \end{equation*}
  and moreover, since $G(x_{n})\to \infty$, we conclude that
  \begin{equation}
    \label{lim:Psqrtn}
    \limsup_{n\to \infty}\frac{\sum_{i=0}^{n-1}P_{i+1}}{\sqrt{n}}=
    \infty.
  \end{equation}
  By applying Lemma \ref{lem:Pmug}, we conclude that $P_{n+1}$ is positive for
  big enough $n$.  Thus $\sum_{i=0}^{n-1}P_{i+1}$ has a limit as $n\to
  \infty$, finite or infinite.  Consequently, (\ref{lim:Psqrtn}) implies that
  the limit is infinite.  Together with Lemma \ref{lem:Pmug} this implies that
  $\sum_{i=0}^\infty |x_i|^{\mu}$ can not be finite on a set of nonzero probability,
  i.e. $\sum_{i=0}^\infty |x_i|^{\mu}=\infty$ a.s.  Therefore we can apply
  Toeplitz Lemma, or Lemma \ref{lem:toepl_cor1}, and obtain that
  \begin{equation*}
  \lim_{n\to \infty}\frac{\sum_{i=0}^{n-1}P_{i+1}}
  {a(\mu+1)\sum_{i=0}^{n-1}|x_i|^{\mu}}
  = \lim_{n\to\infty}\frac{P_{n+1}}{a(\mu+1)|x_n|^{\mu}}=S>0.
  \end{equation*}
  Combining this with (\ref{lim:Psqrtn}) gives
  \begin{equation*}
    \limsup_{n\to \infty}\frac{\sum_{i=0}^{n-1}|x_i|^{\mu}}{\sqrt{n}}=
    \infty.
  \end{equation*}
  Since
  \begin{equation*}
  \sqrt{n}>\frac 12\sum_{i=1}^{n}i^{-\frac 12},
  \end{equation*}
  we can estimate
  \begin{equation*}
  \limsup_{n\to\infty} \frac{\sum_{i=0}^{n-1}|x_i|^{\mu}}
  {\frac12\sum_{i=1}^{n-1}i^{-\frac12}}
  = \limsup_{n\to\infty} \frac{\sum_{i=0}^{n-1}|x_i|^{\mu}}
  {\frac 12\sum_{i=1}^{n}i^{-\frac 12}}
  \ge \limsup_{n\to\infty}
  \frac{\sum_{i=0}^{n-1}|x_i|^{\mu}}{\sqrt{n}} = \infty.
  \end{equation*}
  Applying Lemma \ref{lem:toepl_cor1} again, from the last limit we conclude that
  \begin{equation*}
  \limsup_{n\to\infty}\sqrt{|x_n|^{\mu_g}n}
  = \limsup_{n\to\infty}|x_n|^{\mu}\sqrt{n}=\infty.
  \end{equation*}
\end{proof}

%%%%%%%%%%%%%%%%%%%%%%%%%%%%%%%%%%%%%%

\begin{proof}[Proof of (\ref{eq:nexct2})]
  In the previous part we proved
  that  $\sum_{i=0}^\infty P_{i+1}=\infty$, therefore
  \begin{equation*}
  \liminf_{n\to\infty}\frac{\sum_{i=0}^{n-1}P_{i+1}}{\sqrt{n}}\ge 0.
  \end{equation*}
  After dividing both parts of the decomposition
  (\ref{eq:Gsum}) by $\sqrt{n}$, taking the limsup of both parts and
  applying Lemma \ref{lem:M} we obtain
  \begin{equation*}
    \begin{split}
      \limsup_{n\to\infty}\frac{G(x_{n})}{\sqrt{n}}
      &= \limsup_{n\to\infty}\left[\frac{G(x_0)}{\sqrt{n}}+
        \frac1{\sqrt{n}} \sum_{i=0}^{n-1}P_{i+1}
        + \frac{M_n+T_n}{\sqrt{n}}\right]\\
      &\ge \limsup_{n\to\infty}\frac {M_n+T_n}{\sqrt{n}}=\infty.
    \end{split}
  \end{equation*}
  Therefore,
  \begin{equation*}
    \infty = \limsup_{n\to\infty}\frac{G(x_{n})}{\sqrt{n}}
    = \frac 1{a\mu}\limsup_{n\to\infty}\frac{|x_{n}|^{-\mu}}{\sqrt{n}},
  \end{equation*}
  or
  \begin{equation*}
    \liminf_{n\to\infty}\sqrt{n|x_{n}|^{\mu_g}}=0.
  \end{equation*}
  The theorem is proved.
\end{proof}

%%%%%%%%%%%%%%%%%%%%%%%%%%% Bibliography %%%%%%%%%%%%%%%%%%%%%%%%%%

\end{document}